\documentclass[twocolumn,twoside]{IEEEtran}
\usepackage{amsfonts}
\usepackage{amssymb}
\usepackage{amsmath}
\usepackage{algorithmic}
\usepackage{graphicx}
\usepackage{grffile}

\ifCLASSINFOpdf
\else
\fi
\newtheorem{theorem}{Theorem}
\newtheorem{proposition}[theorem]{Proposition}
\input{tcilatex}
\begin{document}

\title{BM3D frames and variational image deblurring}
\author{Aram~Danielyan, Vladimir~Katkovnik, and Karen~Egiazarian, %
\IEEEmembership{Senior Member,~IEEE} \thanks{%
All authors are with Department of Signal Processing, Tampere University of
Technology, P. O. Box 553, 33101 Tampere, Finland (e-mail:
firstname.lastname@tut.fi).} \thanks{%
This work is supported by Academy of Finland: project no. 213462, 2006-2011
(Finnish Programme for Centres of Excellence in Research) and project no.
138207, 2011-2014, and by Tampere Doctoral Programme in Information Science
and Engineering (TISE).} }
\maketitle

\begin{abstract}
A family of the Block Matching 3-D (BM3D) algorithms for various imaging
problems has been recently proposed within the framework of nonlocal
patch-wise image modeling \cite{BM3D_TIP}, \cite{DEBBM3D_SPIE2008}. In this
paper we construct analysis and synthesis frames, formalizing the BM3D image
modeling and use these frames to develop novel iterative deblurring
algorithms. We consider two different formulations of the deblurring
problem: one given by minimization of the single objective function and
another based on the Nash equilibrium balance of two objective functions.
The latter results in an algorithm where the denoising and deblurring
operations are decoupled. The convergence of the developed algorithms is
proved. Simulation experiments show that the decoupled algorithm derived
from the Nash equilibrium formulation demonstrates the best numerical and
visual results and shows superiority with respect to the state of the art in
the field, confirming a valuable potential of BM3D-frames as an advanced
image modeling tool.
\end{abstract}

\section{Introduction}

\IEEEPARstart{W}{e} consider image restoration from a blurry and noisy
observation. Assuming a circular shift-invariant blur operator and additive
zero-mean white Gaussian noise the conventional observation model is
expressed as%
\begin{equation}
\mathbf{z}=\mathbf{Ay}+\sigma \mathbf{\varepsilon ,}
\label{observations_vect}
\end{equation}%
where $\mathbf{z},\mathbf{y}\in 
%TCIMACRO{\U{211d} }%
%BeginExpansion
\mathbb{R}
%EndExpansion
^{N}$ are vectors representing the observed and true image, respectively, $%
\mathbf{A}$ is an $N\times N$ blur matrix, $\mathbf{\varepsilon }\sim 
\mathcal{N}(\mathbf{0}_{N\times 1},\mathbf{I}_{N\times N})$ is a vector of
i.i.d. Gaussian random variables, and $\sigma $ is the standard deviation of
the noise. The deblurring problem is to reconstruct $\mathbf{y}$ from the
observation $\mathbf{z}$. The most popular approach is to formulate
reconstruction as a variational optimization problem, where the desired
solution minimizes a criterion composed of fidelity and penalty terms. The
fidelity ensures that the solution agrees with the observation, while the
penalty provides regularization of the optimization problem through a prior
image model. Typically, the fidelity term is derived from the negative
log-likelihood function. For the Gaussian observation model (\ref%
{observations_vect}) the fidelity term has the form $\dfrac{1}{2\sigma ^{2}}%
\left\Vert \mathbf{z}-\mathbf{Ay}\right\Vert _{2}^{2}$, and the minimization
criterion is given as%
\begin{equation}
J=\frac{1}{2\sigma ^{2}}\left\Vert \mathbf{z}-\mathbf{Ay}\right\Vert
_{2}^{2}+\tau \cdot pen\left( \mathbf{y}\right) ,  \label{111}
\end{equation}%
where $||\mathbf{\cdot ||}_{2}$ stands for the Euclidean norm, $pen(\cdot )$
is a penalty functional and $\tau >0$ is a regularization parameter.

Image modeling\ lies at the core of image reconstruction problems. Recent
trends are concentrated on \emph{sparse representation} techniques, where
the image is assumed to be defined as a combination of few \emph{atomic}
functions taken from a certain \emph{dictionary}. It follows that the image
can be parameterized and approximated locally or nonlocally by these
functions. To enable sparse approximations, the dictionary should be rich
enough to grasp all variety of the images. Clearly, bases are too limited
for this task and one needs to consider overcomplete systems with a number
of elements essentially larger than the dimensionality of the approximated
images. \emph{Frames} are generalization of the concept of basis to the case
when the atomic functions are linearly dependent and form an overcomplete
system \cite{ChristensenFrames}. There is a vast amount of literature
devoted to the sparsity based models and methods for imaging. An excellent
introduction and overview of this area can be found in the recent book \cite%
{Elad-BOOK}.

The contribution of this paper concerns three main aspects of image
deblurring: image modeling, variational problem formulation, and algorithmic
reconstruction.

First, the BM3D image modeling developed in \cite{BM3D_TIP} is formalized in
terms of the overcomplete sparse frame representation. We construct analysis
and synthesis BM3D-frames and study their properties. The analysis and
synthesis developed in BM3D are interpreted as a general sparse image
modeling applicable to variational formulations of various image processing
problems.

Second, we consider two different formulations of the image deblurring
problem: one given by minimization of the objective function and another
based on the Nash equilibrium. The latter approach results in an algorithm
where the denoising and the deblurring operations are decoupled.

Third, it is shown by simulation experiments that the best image
reconstruction both visually and numerically is obtained by the algorithm
based on decoupling of blur inverse and noise filtering. To the best of our
knowledge, this algorithm provides results which are the state-of-art in the
field.

Here we extend and develop our preliminary ideas sketched in \cite%
{WITMSE-2010}. The BM3D frames are now constructed explicitly, taking into
account the particular form of the 3D transform. Proofs of the frame
properties are presented. We develop algorithms for the analysis and
synthesis-based problem formulations introduced in \cite{WITMSE-2010} and
provide their convergence analysis. The problem formulation based on the
Nash equilibrium and the corresponding decoupled deblurring algorithm are
novel developments.

The paper is organized as follows. We start from a presentation of the BM3D
image modeling and introduce BM3D-frames (Section \ref{sec:BM3D_modeling}).
The variational image reconstruction is a subject of Section \ref%
{sec:Variational_deblurring}. The algorithms based on the analysis and
synthesis formulations are derived in this section. The algorithm based on
the Nash equilibrium is presented in Section \ref{sec:Decoupling}.
Convergence results for the proposed algorithms are given in Section \ref%
{sec:Convergence}. Implementation of the algorithms is discussed in Section %
\ref{sec:Implementations}. The experiments and comparison of the algorithms
are given in Section \ref{sec:Experiments}. In Section \ref{sec:Discussion}
we discuss the principal differences of the decoupled formulation compared
to the analysis and synthesis formulations. Concluding remarks are done in
the last section. Proofs of mathematical statements are given in Appendix.

\section{Overcomplete BM3D image modeling}

\label{sec:BM3D_modeling}

BM3D is a nonlocal image modelling technique based on adaptive, high order
groupwise models. Its detailed discussion can be found in \cite%
{IJCV-Katkonik2010}. Below, using the example of the denoising algorithm 
\cite{BM3D_TIP}, we recall the concept of the BM3D modeling. The denoising
algorithm can be split into three steps.

\begin{enumerate}
\item \emph{Analysis}. Similar image blocks are collected in groups. Blocks
in each group are stacked together to form 3-D data arrays, which are
decorrelated using an invertible 3D transform.

\item \emph{Processing}. The obtained 3-D group spectra are filtered by hard
thresholding.

\item \emph{Synthesis.} The filtered spectra are inverted, providing
estimates for each block in the group. These blockwise estimates are
returned to their original positions and the final image\ reconstruction is
calculated as a weighted average of all the obtained blockwise estimates.
\end{enumerate}

The blocking imposes a localization of the image on small pieces where
simpler models may fit the observations. It has been demonstrated that a
higher sparsity of the signal representation and a lower complexity of the
model can be achieved using joint 3D groupwise instead of 2D blockwise
transforms. This joint 3D transform dramatically improves the effectiveness
of image spectrum approximation.

The total number of groupwise spectrum elements is much larger than the
image size, and we arrive to an \emph{overcomplete} or \emph{redundant} data
approximation. This redundancy is important for effectiveness of the BM3D
modeling.

Our target is to give a strict frame interpretation of the analysis and
synthesis operations in BM3D.

\subsection{Matrix representation of analysis and synthesis operations}

Let $\mathbf{Y}$ be a $\sqrt{N}\times \sqrt{N}$ square matrix representing
the image data and $\mathbf{y}$ be the corresponding $%
%TCIMACRO{\U{211d} }%
%BeginExpansion
\mathbb{R}
%EndExpansion
^{N}$-vector built from the columns of $\mathbf{Y}$. To each $\sqrt{N_{bl}}%
\times \sqrt{N_{bl}}$ square image block we assign unique index equal to the
index of its upper-left corner element (pixel) in $\mathbf{y}$. We denote a
vector of elements of $j$-th block $\mathbf{Y}_{j}$ by $\mathbf{y}_{j}$ and
define $\mathbf{P}_{j}$ as an $N_{bl}\times N$ matrix of indicators $\left[
0,1\right] $ showing which elements of $\mathbf{y}$ belong to the $j$-th
block, so that $\mathbf{y}_{j}=\mathbf{P}_{j}\mathbf{y}$. For the sake of a
notation simplicity, we assume that the number of blocks in each group is
fixed and equal to $K$. Let $J_{r}=\{j_{r,1},...,j_{r,K}\}$ be the set of
indices of the blocks in the $r$-th group, then grouping is completely
defined by the set $J=\left\{ J_{r}:r=1,...,R\right\} $, where $R$ is a
total number of the groups. It is assumed that for each pixel there is at
least one block containing the pixel and entering in some group.

The particular form of the 3-D decorrelating transform constitutes an
important part of the BM3D modeling. It is constructed as a separable
combination of 2-D intrablock and 1-D interblock transforms. The 2-D
transform, in turn, is typically implemented as a separable combination of
1-D transforms. Let $\mathbf{D}_{2}$ and $\mathbf{D}_{1}$ be $\sqrt{N_{bl}}%
\times \sqrt{N_{bl}}$ and $K\times K$ size matrices representing
respectively 1-D interblock and 1-D intrablock transforms. Then the
separable 2-D transform for the block $\mathbf{Y}_{j}$ is given by the
formula%
\begin{equation*}
\mathbf{\Theta }_{j}=\mathbf{D}_{2}\mathbf{Y}_{j}\mathbf{D}_{2}^{T}.
\end{equation*}%
\ \ The vectorization of this formula using the Kronecker matrix product $%
\mathbf{\otimes }$ gives%
\begin{equation*}
\mathbf{\theta }_{j}=\left( \mathbf{D}_{2}\mathbf{\otimes D}_{2}\right)
\cdot \mathbf{y}_{j},
\end{equation*}%
where $\mathbf{\theta }_{j},\mathbf{y}_{j}\in 
%TCIMACRO{\U{211d} }%
%BeginExpansion
\mathbb{R}
%EndExpansion
^{N_{bl}}$ are the vectors corresponding to the matrices $\mathbf{\Theta }%
_{j}$ and $\mathbf{Y}_{j}$, respectively. To obtain the 3-D spectrum of the $%
r$-th group we form the $N_{bl}\times K$ matrix of the vectorized spectrums $%
\left[ \mathbf{\theta }_{j_{r,1}},\mathbf{\theta }_{j_{r,2}},...,\mathbf{%
\theta }_{j_{r,K}}\right] $ and apply the 1-D interblock transform to each
row of this matrix%
\begin{equation*}
\mathbf{\Omega }_{r}=\left[ \mathbf{\theta }_{j_{r,1}},\mathbf{\theta }%
_{j_{r,2}},...,\mathbf{\theta }_{j_{r,K}}\right] \cdot \mathbf{D}_{1}^{T}.
\end{equation*}%
Performing vectorization again, we express the 3-D group spectrum
coefficients in a compact form:%
\begin{eqnarray*}
\mathbf{\omega }_{r} &=&\sum\nolimits_{j\in J_{r}}\mathbf{d}_{j}\otimes %
\left[ \left( \mathbf{D}_{2}\mathbf{\otimes D}_{2}\right) \cdot \mathbf{y}%
_{j}\right] \\
&=&\left( \sum\nolimits_{j\in J_{r}}\mathbf{d}_{j}\otimes \left[ \left( 
\mathbf{D}_{2}\mathbf{\otimes D}_{2}\right) \mathbf{P}_{j}\right] \right)
\cdot \mathbf{y},
\end{eqnarray*}%
where $\mathbf{\omega }_{r}$ is the columnwise vectorized matrix $\mathbf{%
\Omega }_{r}$ and $\mathbf{d}_{j}$ is the $j$-th column of $\mathbf{D}_{1}$.
Finally, denoting%
\begin{equation}
\mathbf{\Phi }_{r}=\sum\nolimits_{j\in J_{r}}\mathbf{d}_{j}\otimes \left[
\left( \mathbf{D}_{2}\mathbf{\otimes D}_{2}\right) \mathbf{P}_{j}\right] ,
\label{spectr1}
\end{equation}%
we express the joint 3D \emph{groupwise spectrum} $\mathbf{\omega }=\left[ 
\mathbf{\omega }_{1}^{T},\ldots ,\mathbf{\omega }_{R}^{T}\right] ^{T}\in 
\mathbb{R}^{M}$ of the image $\mathbf{Y}$ in the vector-matrix form%
\begin{equation}
\mathbf{\omega }=\left[ 
\begin{array}{c}
\mathbf{\Phi }_{1} \\ 
\vdots \\ 
\mathbf{\Phi }_{R}%
\end{array}%
\right] \cdot \mathbf{y}=\mathbf{\Phi y.}  \label{eq:def_PHI}
\end{equation}%
The matrix $\mathbf{\Phi }$ defined by the formulas (\ref{spectr1})-(\ref%
{eq:def_PHI}) gives an explicit representation of the BM3D analysis
operation.

The synthesis matrix is derived similarly. First, the inverse 3-D transform
is applied to each group spectrum $\mathbf{\omega }_{r}$ and then obtained
block estimates are returned to their original positions by $\mathbf{P}%
_{j}^{T},j\in J_{r}$. The estimate obtained from the $r$-th group spectrum
is expressed as $\mathbf{\Psi }_{r}\mathbf{\omega }_{r}$, where%
\begin{equation}
\mathbf{\Psi }_{r}=\sum_{j\in J_{r}}\mathbf{d}_{j}^{T}\otimes \left[ \mathbf{%
P}_{j}^{T}\left( \mathbf{D}_{2}\mathbf{\otimes D}_{2}\right) ^{T}\right]
\label{eq:def_PSI_r}
\end{equation}%
is an $N\times N_{bl}$ matrix.

The final image estimate is defined as the weighted mean of the groupwise
estimates using weights $g_{r}>0$. Hence the synthesis operation has the form%
\begin{equation}
\mathbf{y}=\mathbf{\Psi \omega =W}^{-1}\cdot \lbrack g_{1}\mathbf{\Psi }%
_{1},\ldots ,g_{R}\mathbf{\Psi }_{R}]\cdot \mathbf{\omega },
\label{eq:def_PSI}
\end{equation}%
where%
\begin{equation}
\mathbf{W=}\sum_{r}g_{r}\sum_{j\in J_{r}}\mathbf{P}_{j}^{T}\mathbf{P}_{j}
\label{eq_W}
\end{equation}%
normalizes the weighted mean. $\mathbf{W}$ is a diagonal matrix, since all
products $\mathbf{P}_{j}^{T}\mathbf{P}_{j}$ are diagonal matrices. The $m$%
-th diagonal element of $\mathbf{P}_{j}^{T}\mathbf{P}_{j}$ is $1$ if the $m$%
-th pixel of $\mathbf{y}$ belongs to the $j$-th block, otherwise it is $0$.
Thus, the $m$-th diagonal elements of the matrix-sum $\sum_{j\in I_{r}}%
\mathbf{P}_{j}^{T}\mathbf{P}_{j}$ indicates the number of blocks in the $r$%
-th group containing $m$-th pixel.

The matrix $\mathbf{\Psi }$ defined by the formulas (\ref{eq:def_PSI_r})-(%
\ref{eq_W}) gives the matrix representation of the BM3D synthesis operation.

\subsection{Frame interpretation}

\begin{proposition}
\label{ref_propFrames}\emph{The following equations hold for the matrices }$%
\mathbf{\Phi }$\emph{\ and }$\mathbf{\Psi }$\emph{\ defined by (\ref%
{eq:def_PHI}) and (\ref{eq:def_PSI}):}%
\begin{eqnarray}
&&\mathbf{\Phi }^{T}\cdot \mathbf{\Phi }=\sum_{r}\sum_{j\in I_{r}}\mathbf{P}%
_{j}^{T}\mathbf{P}_{j}>0,  \label{Prop_1_1} \\
&&\mathbf{\Psi }\cdot \mathbf{\Psi }^{T}=\sum_{r}g_{r}^{2}\sum_{j\in I_{r}}%
\mathbf{P}_{j}^{T}\mathbf{P}_{j}\mathbf{W}^{-2}>0,  \label{Prop_1_1_1} \\
&&\mathbf{\Psi }\cdot \mathbf{\Phi }=\mathbf{I}_{N\times N}.
\label{Prop_1_2}
\end{eqnarray}%
The proof is presented in Appendix \ref{app:FrameProps}.
\end{proposition}

It follows from Proposition \ref{ref_propFrames} that rows of $\mathbf{\Phi }
$ constitute a frame $\left\{ \mathbf{\phi }_{n}\right\} $ in $%
%TCIMACRO{\U{211d} }%
%BeginExpansion
\mathbb{R}
%EndExpansion
^{N}$. Indeed, let us verify the frame inequality. Using the analysis
formula $\mathbf{\omega }=\mathbf{\Phi y}$ we obtain%
\begin{align}
& \sum_{n}\left\vert \left\langle \phi _{n},\mathbf{y}\right\rangle
\right\vert ^{2}=\mathbf{\omega }^{T}\mathbf{\omega }=  \notag \\
& =\mathbf{y}^{T}\mathbf{\Phi }^{T}\mathbf{\Phi y}=\mathbf{y}^{T}\mathbf{%
\cdot }\sum_{r}\sum_{j\in I_{r}}\mathbf{P}_{j}^{T}\mathbf{P}_{j}\mathbf{%
\cdot y}.  \label{eq:frame1}
\end{align}%
If $a$ and $b$ are respectively minimum and maximum values of the diagonal
matrix $\sum_{r}\sum_{j\in I_{r}}\mathbf{P}_{j}^{T}\mathbf{P}_{j}$, then for
any $\mathbf{y}\in 
%TCIMACRO{\U{211d} }%
%BeginExpansion
\mathbb{R}
%EndExpansion
^{N}$ holds the frame inequality%
\begin{equation}
a\cdot \left\Vert \mathbf{y}\right\Vert ^{2}\leq \sum_{n}\left\vert
\left\langle \phi _{n},\mathbf{y}\right\rangle \right\vert ^{2}\leq b\cdot
\left\Vert \mathbf{y}\right\Vert ^{2}.  \label{frame2}
\end{equation}%
The frame $\left\{ \mathbf{\phi }_{n}\right\} $ is not tight because $a\neq
b $. This follows from the fact that the elements on the diagonal of matrix $%
\sum_{r}\sum_{j\in I_{r}}\mathbf{P}_{j}^{T}\mathbf{P}_{j}$ count the number
of blocks containing a given pixel. These values are different for different
pixels, since pixels from the blocks possessing higher similarity to other
blocks participate in a larger number of groups.

Similarly, using (\ref{Prop_1_1_1}) we can show that columns of $\mathbf{%
\Psi }$ constitute a non-tight frame $\left\{ \mathbf{\psi }_{n}\right\} $.
From equation (\ref{Prop_1_2}) it follows that $\left\{ \mathbf{\phi }%
_{n}\right\} $ is dual to $\left\{ \mathbf{\psi }_{n}\right\} $. In general $%
\left\{ \mathbf{\phi }_{n}\right\} $ is an alternative dual and becomes
canonical dual only when all weights $g_{r}$ are equal.

We would like to emphasize that since groups and weights are selected data
adaptively, the constructed frames are also data adaptive.

The presented frame interpretation allows to extend the scope of the BM3D
modeling to the modern variational image reconstruction techniques.

\section{Variational image deblurring}

\label{sec:Variational_deblurring}

The frame based variational image reconstruction problem allows two
different formulations depending on what kind of image modeling, analysis or
synthesis is used \cite{Elad-BOOK}. In the \emph{analysis} formulation the
relation between the image and spectrum variables is given by the analysis
equation $\mathbf{\omega }=\mathbf{\Phi y}$. The problem is formalized as a
constrained optimization:%
\begin{equation}
\left( \mathbf{\hat{\omega}},\mathbf{\hat{y}}\right) =\arg \min_{\mathbf{%
\omega ,y}}\{\frac{1}{2\sigma ^{2}}||\mathbf{z}-\mathbf{Ay||}_{2}^{2}+\tau
\cdot \left\Vert \mathbf{\omega }\right\Vert _{p}|\mathbf{\omega }=\mathbf{%
\Phi y\}},  \label{var_analysis}
\end{equation}%
where $\left\Vert \cdot \right\Vert _{p}$ is the standard notation of the $%
l_{p}$-norm.

In the \emph{synthesis}\textit{\ }formulation the relation is given by the
synthesis equation\ $\mathbf{y=\Psi \omega },$ leading to the constrained
optimization:%
\begin{equation}
\left( \mathbf{\hat{\omega},\hat{y}}\right) =\arg \min_{\mathbf{\omega },%
\mathbf{y}}\{\frac{1}{2\sigma ^{2}}||\mathbf{z}-\mathbf{Ay||}_{2}^{2}+\tau
\cdot \left\Vert \mathbf{\omega }\right\Vert _{p}|\mathbf{y=\Psi \omega \}}.
\label{var_synthesis}
\end{equation}%
These problems have equivalent unconstrained forms in which they usually
encounter in literature. To obtain them it is enough to eliminate $\mathbf{%
\omega }$ and $\mathbf{y}$ respectively from (\ref{var_analysis}) and (\ref%
{var_synthesis}). The analysis problem is then formulated as the
minimization in the image domain%
\begin{equation}
\mathbf{\hat{y}}=\arg \min_{\mathbf{y}}\{\frac{1}{2\sigma ^{2}}|\mathbf{%
|z-Ay||}_{2}^{2}+\tau \cdot \left\Vert \mathbf{\Phi y}\right\Vert _{p}%
\mathbf{\}}.  \label{st_analysis}
\end{equation}%
Similarly, the synthesis problem is formulated as the minimization in the
spectrum domain%
\begin{equation}
\mathbf{\hat{\omega}}=\arg \min_{\mathbf{\omega }}\{\frac{1}{2\sigma ^{2}}||%
\mathbf{z-A\Psi \omega ||}_{2}^{2}+\tau \cdot \left\Vert \mathbf{\mathbf{%
\omega }}\right\Vert _{p}\mathbf{\}}.  \label{st_synthesis}
\end{equation}

Despite of the algebraic similarity, the analysis and synthesis formulations
generally lead to different solutions. A detailed discussion of the
nontrivial connections between the analysis and synthesis formulations can
be found in \cite{Elad-SynthesisAnalysis}.

The problems (\ref{var_analysis})-(\ref{st_synthesis}) and the corresponding
solution techniques recently become a subject of an intensive study. In
particular, several algorithms have been suggested for the convex $l_{1}$%
-norm penalty. These algorithms sharing many common ideas are known under
different names such as \emph{split Bregman iterations}\textit{\ }\cite%
{SIAM-IM-OSHER-2009}\textit{, \emph{iterative shrinkage algorithms }}\cite%
{Beck-Teboulle2009}\textit{, }\emph{alternating direction method of
multipliers} \cite{Afonso-Dias-Splitting}, \emph{majorization-minimization}%
\textit{\ }algorithms \cite{Dias-2009-TVMM}. In this paper similar to \cite%
{Afonso-Dias-AL} we confine ourself to the Augmented Langrangian (AL)
technique, using it as a simple and efficient tool for an explicit
derivation of the reconstruction algorithms. This AL\ technique, introduced
independently by Hestenes \cite{Hestens-Lagrange} and Powell \cite%
{Powell-Lagrange} is now widely used for minimization of convex functionals
under linear equality constraints.

\subsection{Analysis-based reconstruction}

\label{sec:AnalysisBased}

The AL criterion for the \emph{analysis}\textit{\ }formulation (\ref%
{var_analysis})\textit{\ }takes the form:%
\begin{eqnarray}
L_{\text{a}}\left( \mathbf{y,\omega ,\lambda }\right) &=&\frac{1}{2\sigma
^{2}}\left\Vert \mathbf{z}-\mathbf{Ay}\right\Vert _{2}^{2}+\tau \cdot
\left\Vert \mathbf{\omega }\right\Vert _{p}+  \notag \\
&&\frac{1}{2\gamma }\left\Vert \mathbf{\omega }-\mathbf{\Phi y}\right\Vert
_{2}^{2}+\frac{1}{\gamma }\left\langle \mathbf{\omega }-\mathbf{\Phi
y,\lambda }\right\rangle ,  \label{A_AL_1}
\end{eqnarray}%
where $\mathbf{\lambda }$ is a vector of the Lagrange multipliers, $\gamma
>0 $ is a parameter and the subscript 'a' indicates the analysis
formulation. The saddle problem associated with the Lagrangian $L_{\text{a}}$
provides the solution of the constrained optimization problem (\ref%
{var_analysis}).

Finding the saddle point\ requires minimization of $L_{\text{a}}$ with
respect to the variables $\mathbf{y},\mathbf{\omega }$ and maximization with
respect to $\mathbf{\lambda }$. A common practical approach is to find the
saddle point by performing alternating optimization. Applied to (\ref{A_AL_1}%
) it results in the following iterative scheme:

\emph{Repeat for }$t=0,1,...$%
\begin{eqnarray}
\mathbf{y}_{t+1} &=&\arg \min_{\mathbf{y}}L_{\text{a}}\left( \mathbf{y},%
\mathbf{\omega }_{t},\mathbf{\lambda }_{t}\right) ,  \label{AUG_1s} \\
\mathbf{\omega }_{t+1} &=&\arg \min_{\mathbf{\omega }}L_{\text{a}}\left( 
\mathbf{y}_{t+1},\mathbf{\omega },\mathbf{\lambda }_{t}\right) ,
\label{AUG_2s} \\
\mathbf{\lambda }_{t+1} &=&\mathbf{\lambda }_{t}+\beta \cdot \left( \mathbf{%
\omega }_{t+1}-\mathbf{\Phi y}_{t+1}\right) ,  \label{AUG_3s}
\end{eqnarray}

\emph{until convergence.}

\noindent Here maximization with respect to $\mathbf{\lambda }$ is produced
as a step (\ref{AUG_3s}) in the direction of the gradient $\nabla _{\mathbf{%
\lambda }}L_{\text{a}}$, with a step-size $\beta >0$. The convergence of the
scheme (\ref{AUG_1s})-(\ref{AUG_3s}) is studied in \cite{SIAM-IM-OSHER-2009}.

\emph{Minimization with respect to }$\mathbf{y}$. Since $L_{\text{a}}$ is
quadratic with respect to $\mathbf{y}$ the optimal solution is defined by
the linear equation%
\begin{equation}
\left( \frac{1}{\sigma ^{2}}\mathbf{A}^{T}\mathbf{A}+\frac{1}{\gamma }%
\mathbf{\Phi }^{T}\mathbf{\Phi }\right) \cdot \mathbf{y=}\frac{1}{\sigma ^{2}%
}\mathbf{\mathbf{A}^{T}\mathbf{z}}+\frac{1}{\gamma }\mathbf{\Phi }^{T}\left( 
\mathbf{\omega }+\mathbf{\lambda }\right) .  \label{A-AL_3}
\end{equation}%
\ \ We denote by $\hat{Y}_{\text{a}}\left( \mathbf{\omega ,\mathbf{\lambda }}%
\right) $ the operator giving the solution of (\ref{A-AL_3}).

\emph{Minimization with respect to }$\mathbf{\omega }$\textbf{. }Regrouping
the terms in $L_{\text{a}}$ we arrive to the following formula%
\begin{eqnarray*}
L_{\text{a}}\left( \mathbf{y,\omega ,\lambda }\right) &=&\frac{1}{2\sigma
^{2}}\left\Vert \mathbf{z}-\mathbf{Ay}\right\Vert _{2}^{2}+\tau \cdot
\left\Vert \mathbf{\omega }\right\Vert _{p}+ \\
&&\frac{1}{2\gamma }\left\Vert \mathbf{\omega }-\left( \mathbf{\Phi y-%
\mathbf{\lambda }}\right) \right\Vert _{2}^{2}-\frac{1}{2\gamma }\left\Vert 
\mathbf{\lambda }\right\Vert _{2}^{2}.
\end{eqnarray*}%
\ \ Since the first and the last terms do not depend on $\mathbf{\omega }$,
the problem is reduced to the optimization%
\begin{equation}
\mathbf{\hat{\omega}}=\arg \min_{\mathbf{\omega }}\tau \cdot \left\Vert 
\mathbf{\omega }\right\Vert _{p}+\frac{1}{2\gamma }\left\Vert \mathbf{\omega 
}-\left( \mathbf{\Phi y-\mathbf{\lambda }}\right) \right\Vert _{2}^{2}.
\label{op_A_AL_omega}
\end{equation}%
\ For $p\leq 1,$ the $l_{p}$-norm is non-differentiable which makes
optimization on $\mathbf{\omega }$ non-trivial. Nevertheless, for $p=0$ and $%
p=1$ there are well known analytical solutions.

Let us denote $\mathbf{\mathbf{b}}=\mathbf{\Phi y-\mathbf{\lambda }}$, then (%
\ref{op_A_AL_omega}) takes the form 
\begin{equation}
\mathbf{\hat{\omega}}=\arg \min_{\mathbf{\omega }}\tau \cdot \left\Vert 
\mathbf{\omega }\right\Vert _{p}+\frac{1}{2}\left\Vert \mathbf{\omega }-%
\mathbf{b}\right\Vert _{2}^{2},\mathbf{\omega },\mathbf{b\in }%
%TCIMACRO{\U{211d} }%
%BeginExpansion
\mathbb{R}
%EndExpansion
^{M}.  \label{eq:prob_thresh}
\end{equation}%
Depending on the used norm the solution of (\ref{eq:prob_thresh}) is given
either by the hard or soft thresholding according to the formula:%
\begin{eqnarray}
\mathbf{\hat{\omega}} &=&\mathfrak{Th}_{\tau }\left( \mathbf{b}\right) = 
\notag \\
&&\hspace{-0.8cm}\left\{ 
\begin{array}{l}
\mathfrak{Th}_{\tau }^{soft}\left( \mathbf{b}\right) =sign\left( \mathbf{b}%
\right) \circ \max \left( \left\vert \mathbf{b}\right\vert -\tau ,0\right)
,~p=1, \\ 
\mathfrak{Th}_{\sqrt{2\tau }}^{hard}\left( \mathbf{b}\right) =\mathbf{b}%
\circ 1\left( \left\vert \mathbf{b}\right\vert \geq \sqrt{2\tau }\right)
,~p=0.%
\end{array}%
\right.  \label{thresholding_func}
\end{eqnarray}%
Here all vector operations are elementwise, and '$\circ $' stands for the
elementwise product of two vectors. We use $\mathfrak{Th}_{\tau }\left( 
\mathbf{b}\right) $ as a generic notation for the thresholding operator.
Note, that for a given $\tau $ the thresholding levels for the hard and soft
thresholdings are calculated differently.

Applying the general formula (\ref{thresholding_func}) to (\ref%
{op_A_AL_omega}) we obtain the solution in the form 
\begin{equation}
\mathbf{\hat{\omega}}=\mathfrak{Th}_{\tau \gamma }\left( \mathbf{\Phi y}-%
\mathbf{\lambda }\right) .  \label{A-AL_2}
\end{equation}

Following (\ref{AUG_1s})-(\ref{AUG_3s}) and using (\ref{A-AL_3}) and (\ref%
{A-AL_2}) we define the analysis-based iterative algorithm which is
presented in Figure \ref{alg_Analysis}. In each iteration it first updates
the image estimate using the linear filtering (\ref{A-AL_3}). Then, the
difference between the spectrum $\mathbf{\Phi y}_{t}$ and $\mathbf{\lambda }%
_{t}$ is thresholded, what corresponds to the optimization with respect to $%
\mathbf{\omega }$. Finally, the Lagrange multipliers are updated in the
direction of the gradient $\mathbf{\omega }_{t+1}-\mathbf{\Phi y}_{t+1}$.
Process is iterated until some convergence criteria is satisfied.
Particularly, the iterations can be stopped as soon as the difference
between consecutive estimates becomes small enough.

%TCIMACRO{\TeXButton{bF}{\begin{figure}}}%
%BeginExpansion
\begin{figure}%
%EndExpansion
%TCIMACRO{\TeXButton{bAlg}{\begin{algorithmic}}}%
%BeginExpansion
\begin{algorithmic}%
%EndExpansion

\STATE\textbf{input:} $\mathbf{z},\mathbf{A},\mathbf{y}_{\text{init}}$

\STATE\textbf{initialization:}

\STATE\emph{using} $\mathbf{y}_{\text{init}}$ \emph{construct operators} $%
\mathbf{\Phi }$ \emph{and} $\mathbf{\Phi }^{T}$

\STATE\emph{set:} $\mathbf{y}_{0},\mathbf{\omega }_{0},\mathbf{\lambda }_{0}$

\STATE$t=0$

\REPEAT

\STATE$\mathbf{y}_{t+1}\mathbf{=}\hat{Y}_{\text{a}}\left( \mathbf{\omega }%
_{t}\mathbf{,\mathbf{\lambda }}_{t}\right) $

\STATE$\mathbf{\omega }_{t+1}=\mathfrak{Th}_{\tau \gamma }\left( \mathbf{%
\Phi y}_{t+1}-\mathbf{\lambda }_{t}\right) $

\STATE$\mathbf{\lambda }_{t+1}=\mathbf{\lambda }_{t}+\beta \cdot \left( 
\mathbf{\omega }_{t+1}-\mathbf{\Phi y}_{t+1}\right) $

\STATE$t=t+1$

\UNTIL {convergence.}

%TCIMACRO{\TeXButton{eAlg}{\end{algorithmic}}}%
%BeginExpansion
\end{algorithmic}%
%EndExpansion
%TCIMACRO{%
%\TeXButton{Analysis based algorithm}{\caption{Analysis-based deblurring algorithm}}}%
%BeginExpansion
\caption{Analysis-based deblurring algorithm}%
%EndExpansion
\label{alg_Analysis}%
%TCIMACRO{\TeXButton{eF}{\end{figure}}}%
%BeginExpansion
\end{figure}%
%EndExpansion

\subsection{Synthesis-based reconstruction}

\label{sec:SynthesisBased}

The AL criterion for the \emph{synthesis}\textit{\ }formulation (\ref%
{var_synthesis})\textit{\ }takes form:%
\begin{eqnarray}
L_{\text{s}}\left( \mathbf{y},\mathbf{\omega },\mathbf{\lambda }\right) &=&%
\frac{1}{2\sigma ^{2}}\left\Vert \mathbf{z}-\mathbf{Ay}\right\Vert
_{2}^{2}+\tau \cdot \left\Vert \mathbf{\omega }\right\Vert _{p}+  \notag \\
&&\hspace{-0.5cm}\frac{1}{2\gamma }\left\Vert \mathbf{y}-\mathbf{\Psi \omega 
}\right\Vert _{2}^{2}+\frac{1}{\gamma }\left\langle \mathbf{y}-\mathbf{\Psi
\omega ,\lambda }\right\rangle .  \label{S_AL_1}
\end{eqnarray}%
In $L_{\text{s}}$, as opposed to $L_{\text{a}}$, the spectrum variable $%
\mathbf{\omega }$ enters the quadratic term with a matrix factor $\mathbf{%
\Psi }$. It makes the thresholding formula (\ref{thresholding_func})
inapplicable for minimizing $L_{\text{s}}\left( \mathbf{y},\mathbf{\omega },%
\mathbf{\lambda }\right) $ with respect to $\mathbf{\omega }$. One option is
to apply one of the \emph{iterative shrinkage} methods \cite{Elad-BOOK}, but
we prefer to follow a different approach which leads to a simpler solution.
We modify (\ref{S_AL_1}) by introducing a splitting variable $\mathbf{u}\in
R^{M},$ used as an auxiliary estimate of the spectrum $\omega $. The
modified AL takes the form:%
\begin{eqnarray}
\tilde{L}_{\text{s}}\left( \mathbf{y,\omega ,\lambda ,u}\right) &=&\frac{1}{%
2\sigma ^{2}}\left\Vert \mathbf{z}-\mathbf{Ay}\right\Vert _{2}^{2}+\tau
\cdot \left\Vert \mathbf{\omega }\right\Vert _{p}+  \notag \\
&&\frac{1}{2\gamma }\left\Vert \mathbf{y-\Psi u}\right\Vert _{2}^{2}+\frac{1%
}{\gamma }\left\langle \mathbf{y-\Psi u,\lambda }\right\rangle +  \notag \\
&&\frac{1}{2\xi }\left\Vert \mathbf{\omega -u}\right\Vert _{2}^{2}.
\label{S_AL_11}
\end{eqnarray}%
The corresponding saddle point problem is%
\begin{equation}
\arg \min_{\mathbf{y,\omega ,u}}\max_{\mathbf{\lambda }}\tilde{L}_{\text{s}%
}\left( \mathbf{y,\omega ,\lambda ,u}\right) ,  \label{eq_S_AL_2}
\end{equation}%
where optimization with respect to the splitting variable $\mathbf{u}$ is
required.

With a small enough $\xi >0$ penalization by $\dfrac{1}{2\xi }\left\Vert 
\mathbf{\omega -u}\right\Vert _{2}^{2}$ results in $\left\Vert \mathbf{%
\omega -u}\right\Vert _{2}^{2}\rightarrow 0$ what makes the problem (\ref%
{eq_S_AL_2}) equivalent to the saddle problem for (\ref{S_AL_1}). As in the
analysis case we seek for the solution of (\ref{eq_S_AL_2}) by the
alternating optimization of $\tilde{L}_{\text{s}}\left( \mathbf{y,\omega
,\lambda ,u}\right) $ with respect to the variables $\mathbf{y,\omega ,u}$
and $\mathbf{\lambda }$.

\emph{Minimization with respect to}\textit{\ }$\mathbf{y}$ is given by the
solution of the linear equation%
\begin{equation}
\left( \frac{1}{\sigma ^{2}}\mathbf{A}^{T}\mathbf{A}+\frac{1}{\gamma }%
\mathbf{I}_{N\times N}\right) \cdot \mathbf{y=}\frac{1}{\sigma ^{2}}\mathbf{%
\mathbf{A}}^{T}\mathbf{\mathbf{z}}+\frac{1}{\gamma }\left( \mathbf{\Psi u}-%
\mathbf{\lambda }\right) .  \label{S_AL_3}
\end{equation}

\emph{Minimization with respect to }$\mathbf{u}$ satisfies the linear
equation%
\begin{equation}
\left( \frac{1}{\gamma }\mathbf{\Psi }^{T}\mathbf{\Psi }+\frac{1}{\xi }%
\mathbf{I}_{M\times M}\right) \cdot \mathbf{u}=\frac{1}{\gamma }\mathbf{\Psi 
}^{T}\left( \mathbf{y+\mathbf{\lambda }}\right) +\frac{1}{\xi }\mathbf{%
\omega }.  \label{S_AL_3_2}
\end{equation}

\emph{Minimization with respect to}\textit{\ }$\mathbf{\omega }$, thanks to
the splitting variable $\mathbf{u}$, can be obtained by the thresholding (%
\ref{thresholding_func}) with the parameter $\tau \xi $:%
\begin{equation}
\mathbf{\hat{\omega}}=\mathfrak{Th}_{\tau \xi }\left( \mathbf{u}\right) .
\label{S_AL_1_1}
\end{equation}%
We denote by $\hat{Y}_{\text{s}}\left( \mathbf{u,\omega ,\lambda }\right) $
and $\hat{U}_{\text{s}}\mathbf{(y,}\mathbf{\omega },\mathbf{\mathbf{\lambda }%
)}$ the operators giving the solutions of (\ref{S_AL_3}) and (\ref{S_AL_3_2}%
).

Using (\ref{S_AL_3})-(\ref{S_AL_1_1}) we define the synthesis-based
iterative deblurring algorithm which is presented in Figure \ref%
{alg_Synthesis}. At the first two steps the estimates for the image $\mathbf{%
y}_{t}$ and the splitting variable $\mathbf{u}_{t}$ are updated by solving (%
\ref{S_AL_3}) and (\ref{S_AL_3_2}). Then, the splitting variable $\mathbf{u}%
_{t+1}$ is thresholded reducing the complexity of the spectrum estimate $%
\mathbf{\omega }$. Finally, the Lagrange multipliers are updated in the
direction of the gradient $\mathbf{y}_{t+1}\mathbf{-\Psi u}_{t+1}$. Process
is iterated until some convergence criteria is satisfied.

%TCIMACRO{\TeXButton{bF}{\begin{figure}}}%
%BeginExpansion
\begin{figure}%
%EndExpansion
%TCIMACRO{\TeXButton{bAlg}{\begin{algorithmic}}}%
%BeginExpansion
\begin{algorithmic}%
%EndExpansion

\STATE\textbf{input:} $\mathbf{z},\mathbf{A},\mathbf{y}_{\text{init}}$

\STATE\textbf{initialization:}

\STATE\emph{using} $\mathbf{y}_{\text{init}}$ \emph{construct operators} $%
\mathbf{\Psi }$ \emph{and} $\mathbf{\Psi }^{T}$

\STATE\emph{set:} $\mathbf{y}_{0},\mathbf{\omega }_{0},\mathbf{\lambda }_{0},%
\mathbf{u}_{0}$

\STATE$t=0$

\REPEAT

\STATE$\mathbf{y}_{t+1}\mathbf{=}\hat{Y}_{\text{s}}\left( \mathbf{u}_{t}%
\mathbf{,\omega }_{t,}\mathbf{\mathbf{\lambda }}_{t}\right) $

\STATE$\mathbf{u}_{t+1}\mathbf{=}\hat{U}_{\text{s}}\left( \mathbf{y}_{t}%
\mathbf{,\omega }_{t+1},\mathbf{\mathbf{\lambda }}_{t}\right) $

\STATE$\mathbf{\omega }_{t+1}\mathbf{=}\mathfrak{Th}_{\tau \xi }\left( 
\mathbf{u}_{t+1}\right) $

\STATE$\mathbf{\lambda }_{t+1}=\mathbf{\lambda }_{t}+\beta \cdot \left( 
\mathbf{y}_{t+1}\mathbf{-\Psi u}_{t+1}\right) $

\STATE$t=t+1$

\UNTIL {convergence.}

%TCIMACRO{\TeXButton{eAlg}{\end{algorithmic}}}%
%BeginExpansion
\end{algorithmic}%
%EndExpansion
%TCIMACRO{%
%\TeXButton{Synthesis based algorithm}{\caption{Synthesis-based deblurring algorithm}}}%
%BeginExpansion
\caption{Synthesis-based deblurring algorithm}%
%EndExpansion
\label{alg_Synthesis}%
%TCIMACRO{\TeXButton{eF}{\end{figure}}}%
%BeginExpansion
\end{figure}%
%EndExpansion

\section{Decoupling of blur inversion and denoising}

\label{sec:Decoupling}

Above we considered algorithms based on the minimization of a single
objective function. In this section we present an alternative approach based
on formulation of the deblurring as a Nash equlibrium problem for two
objective functions. This approach allows to split the deblurring problem
into two subproblems: a blur inversion and denoising, which are then solved
sequentially. Such a decoupling has several advantages:

\begin{enumerate}
\item The decoupled algorithms are simpler in design and parameter selection;

\item The blur inversion can be implemented efficiently using Fast Fourier
Transform (FFT);

\item Various denoising algorithms can be used in this scheme selected
independently with respect to deblurring;

\item In many cases decoupled algorithms demonstrate better performance than
the algorithms where deblurring and denoising are performed jointly.
\end{enumerate}

Examples of the decoupled deblurring can be found in works \cite%
{DEBBM3D_SPIE2008}, \cite{FORWARD2004TIP}, \cite{IMDEBLUR_SADCT06} and \cite%
{PORTILLA_2007}, where the regularized inverse is followed by different
types of filtering (wavelet, shape-adaptive DCT, BM3D, pyramidal). An
interesting development of this technique is demonstrated in \cite%
{Decouple_China2008} where an iterative algorithm is derived by alternating
optimization of multiple objective functions.

\subsection{Deblurring as a Nash equilibrium problem}

Let us formulate the deblurring problem as the following constrained
optimization:%
\begin{equation}
\left\{ 
\begin{array}{l}
\mathbf{y}^{\ast }=\arg \min\limits_{\mathbf{y}}\frac{1}{2\sigma ^{2}}%
\left\Vert \mathbf{z}-\mathbf{Ay}\right\Vert _{2}^{2}\text{ \hspace{-0.1cm}%
subject \hspace{-0.1cm}to \hspace{-0.1cm}}\left\Vert \mathbf{y-\Psi \omega }%
^{\ast }\right\Vert _{2}^{2}\leq \text{\hspace{-0.1cm}}\varepsilon _{1}, \\ 
\mathbf{\omega }^{\ast }=\arg \min\limits_{\mathbf{\omega }}\tau \cdot
\left\Vert \mathbf{\omega }\right\Vert _{p}\text{ subject to }\left\Vert 
\mathbf{\omega }-\mathbf{\Phi y}^{\ast }\right\Vert _{2}^{2}\leq \varepsilon
_{2},%
\end{array}%
\right.  \label{eq:DecoupledConstrained}
\end{equation}%
where $\varepsilon _{1},\varepsilon _{2}>0$. This problem can be replaced by
the equivalent unconstrained one:%
\begin{equation}
\left\{ 
\begin{array}{c}
\mathbf{y}^{\ast }=\arg \min\limits_{\mathbf{y}}L_{\text{inv}}\left( \mathbf{%
y},\mathbf{\omega }^{\ast }\right) \\ 
\mathbf{\omega }^{\ast }=\arg \min\limits_{\mathbf{\omega }}L_{\text{den}%
}\left( \mathbf{y}^{\ast },\mathbf{\omega }\right)%
\end{array}%
\right. ,  \label{var_decoupled}
\end{equation}%
where%
\begin{eqnarray}
L_{\text{inv}}(\mathbf{y,\omega }) &=&\frac{1}{2\sigma ^{2}}\left\Vert 
\mathbf{z}-\mathbf{Ay}\right\Vert _{2}^{2}+\frac{1}{2\gamma }\left\Vert 
\mathbf{y-\Psi \omega }\right\Vert _{2}^{2},  \label{eq_Decoupled1} \\
L_{\text{den}}(\mathbf{y,\omega }) &=&\tau \cdot \left\Vert \mathbf{\omega }%
\right\Vert _{p}+\frac{1}{2\xi }\left\Vert \mathbf{\omega }-\mathbf{\Phi y}%
\right\Vert _{2}^{2}.  \label{eq_Decoupled2}
\end{eqnarray}%
and $\gamma ,\xi $ are constants selected correspondingly to the values of $%
\varepsilon _{1},\varepsilon _{2}.$

In terms of the game theory the problem (\ref{var_decoupled}) can be
interpreted as a game of two players identified, respectively, with two
variables $\mathbf{y}$ and $\mathbf{\omega }$ \cite%
{Essentials_of_Game_Theory},\cite{GeneralizedNash}. An interaction between
the players is noncooperative because minimization of $L_{\text{inv}}(%
\mathbf{y,\omega })$ with respect to $\mathbf{y}$ in general results in
increase of $L_{\text{den}}(\mathbf{y,\omega })$ and minimization of $L_{%
\text{den}}(\mathbf{y,\omega })$ with respect to $\mathbf{\omega }$
increases $L_{\text{inv}}(\mathbf{y,\omega })$. The equilibrium of this game
called \emph{Nash equilibrium} defines the \emph{fixed point} $\left( 
\mathbf{y}^{\ast },\mathbf{\omega }^{\ast }\right) $ of the optimization.
For $p=1$, problem (\ref{var_decoupled}) is convex.

The objective functions $L_{\text{inv}}$ and $L_{\text{den}}$ allow the
following interpretation. In $L_{\text{inv}}$ the fidelity term $\dfrac{1}{%
2\sigma ^{2}}\left\Vert \mathbf{z}-\mathbf{Ay}\right\Vert _{2}^{2}$
evaluates the divergency between the observation $\mathbf{z}$ and its
prediction $\mathbf{Ay}$. This fidelity is penalized by the norm $\left\Vert 
\mathbf{y-\Psi \omega }\right\Vert _{2}^{2}$ defining a difference between $%
\mathbf{y}$ and its prediction $\mathbf{\Psi \omega }$ through $\mathbf{%
\omega }$. The term $\dfrac{1}{2\xi }\left\Vert \mathbf{\omega }-\mathbf{%
\Phi y}\right\Vert _{2}^{2}$ in $L_{\text{den}}$ evaluates a difference
between the spectrum $\mathbf{\omega }$ and the spectrum prediction $\mathbf{%
\Phi y}$ obtained from $\mathbf{y}$. The error between $\mathbf{\omega }$
and $\mathbf{\Phi y}$ is penalized by the norm $\left\Vert \mathbf{\omega }%
\right\Vert _{p}$.

Hence the Nash equilibrium provides a balance between the fit of the
reconstruction $\mathbf{y}$ to the observation $\mathbf{z}$ and the
complexity of the model $\left\Vert \mathbf{\omega }\right\Vert _{p}$. This
can be contrasted with the analysis and synthesis-based problem formulations
where the balance is provided within a single criterion. As we demonstrate
later the form of the balance plays an essential role in the reconstructions
with non-tight frames.

\subsection{IDD-BM3D algorithm}

To solve (\ref{var_decoupled}) we consider the following iterative procedure:%
\begin{equation}
\left\{ 
\begin{array}{l}
\mathbf{y}_{t+1}=\arg \min\limits_{\mathbf{y}}L_{\text{inv}}(\mathbf{%
y,\omega }_{t}) \\ 
\mathbf{\omega }_{t+1}=\arg \min\limits_{\mathbf{\omega }}L_{\text{den}%
}\left( \mathbf{y}_{t+1},\mathbf{\omega }\right)%
\end{array}%
,\text{ }t=0,1,...\right. .  \label{decoupled}
\end{equation}%
The iterative algorithm (\ref{decoupled}) models the selfish behavior, where
each variable minimizes only its own objective function. These iterations
converge to the fixed point $\left( \mathbf{y}^{\ast },\mathbf{\omega }%
^{\ast }\right) $ of (\ref{var_decoupled}), the corresponding result is
formulated in Section \ref{sec:Convergence}.

Minimization of $L_{\text{inv}}$\ with respect to $\mathbf{y}$ is given by
the solution of the linear equation%
\begin{equation}
\left( \frac{1}{\sigma ^{2}}\mathbf{A}^{T}\mathbf{A}+\frac{1}{\gamma }%
\mathbf{I}\right) \cdot \mathbf{y=}\frac{1}{\sigma ^{2}}\mathbf{\mathbf{A}}%
^{T}\mathbf{\mathbf{z}}+\frac{1}{\gamma }\mathbf{\Psi \omega .}
\label{eq:Dec_AL_1}
\end{equation}%
This step performs regularized inversion of the blur operator.

The minimization of $L_{\text{den}}$\ with respect to $\mathbf{\omega }$ is
obtained by thresholding with the\ threshold parameter $\tau \xi $:%
\begin{equation}
\mathbf{\omega }=\mathfrak{Th}_{\tau \xi }\left( \mathbf{\Phi }\mathbf{y}%
\right) .  \label{eq:Dec_AL_2}
\end{equation}%
Thus, in (\ref{decoupled}) the blur inversion and the denoising steps are
fully decoupled.

The algorithm based on (\ref{decoupled}) is presented in Figure \ref%
{alg_Decoupled}. We call this algorithm Iterative Decoupled Deblurring BM3D
(IDD-BM3D).\footnote{%
We wish to note that IDD-BM3D is similar but not identical to our Augmented
Lagrangian BM3D deblurring (AL-BM3D-DEB) algorithm presented earlier in \cite%
{WITMSE-2010}. The AL-BM3D-DEB algorithm is derived from the analysis-based
formulation (\ref{A_AL_1}). The regularized inverse step (\ref{A-AL_3}) in
AL-BM3D-DEB is replaced by the inverse (\ref{S_AL_3}) obtained from the
synthesis-based formulation (\ref{S_AL_1}). In \cite{WITMSE-2010} this
replacement is treated as an approximation and is not mathematically
rigorous. The presence of the Lagrange multipliers discriminates the
AL-BM3D-DEB algorithm from the IDD-BM3D.}

%TCIMACRO{\TeXButton{bF}{\begin{figure}}}%
%BeginExpansion
\begin{figure}%
%EndExpansion
%TCIMACRO{\TeXButton{bAlg}{\begin{algorithmic}}}%
%BeginExpansion
\begin{algorithmic}%
%EndExpansion

\STATE\textbf{input:} $\mathbf{z},\mathbf{A},\mathbf{y}_{\text{init}}$

\STATE\textbf{initialization:}

\STATE\emph{using} $\mathbf{y}_{\text{init}}$ \emph{construct operators} $%
\mathbf{\Phi }$ \emph{and} $\mathbf{\Psi }$

\STATE\emph{set:} $\mathbf{y}_{0},\mathbf{\omega }_{0}=\mathbf{\Phi y}_{0}$

\STATE$t=0$

\REPEAT

\STATE\emph{Deblurring}\textbf{:}

$\mathbf{y}_{t+1}=\hspace{-0.05cm}\left[ \tfrac{1}{\sigma ^{2}}\mathbf{A}^{T}%
\mathbf{A}+\frac{1}{\gamma }\mathbf{I}\right] ^{-1}\hspace{-0.15cm}\times 
\hspace{-0.05cm}\left[ \frac{1}{\sigma ^{2}}\mathbf{\mathbf{A}}^{T}\mathbf{z}%
+\frac{1}{\gamma }\mathbf{\Psi \omega }_{t}\right] $

\STATE\emph{Denoising}\textbf{:}

$\mathbf{\omega }_{t+1}\mathbf{=}\mathfrak{Th}_{\tau \xi }\left( \mathbf{%
\Phi }\mathbf{y}_{t+1}\right) $

\STATE$t=t+1$

\UNTIL {convergence.}

%TCIMACRO{\TeXButton{eAlg}{\end{algorithmic}}}%
%BeginExpansion
\end{algorithmic}%
%EndExpansion
%TCIMACRO{%
%\TeXButton{IDD-BM3D algorithm}{\caption{IDD-BM3D - Iterative Decoupled Deblurring BM3D algorithm}}}%
%BeginExpansion
\caption{IDD-BM3D - Iterative Decoupled Deblurring BM3D algorithm}%
%EndExpansion
\label{alg_Decoupled}%
%TCIMACRO{\TeXButton{eF}{\end{figure}}}%
%BeginExpansion
\end{figure}%
%EndExpansion

\section{Convergence}

\label{sec:Convergence}

\subsection{Analysis and synthesis-based algorithms}

The main motivation of the AL technique is to replace a constrained
optimization with a simpler saddle-point problem. The equivalence of these%
\textbf{\ }two problems is not a given fact. The classical results stating
equivalence are formulated for the convex and differentiable functions \cite%
{Bertsekas-BOOK}. Since $l_{p}$-norms with $p\leq 1$ are non-differentiable
these results are inapplicable. Nevertheless, for the $l_{1}$-norm the
equivalence can be shown, provided that the constraints in the problem are
linear. In the recent paper \cite{Tai-Wu-2010} the equivalence statement is
proved for the total variation penalty. This proof remains valid for any
convex and non-differentiable penalties, in particularly for the $l_{1}$%
-norm based penalties. The equivalence result is formulated as following:

$(\mathbf{\hat{y},\hat{\omega}})$\emph{\ is a solution of the analysis or
synthesis problems if and only if there exist a saddle-point of the
corresponding ALs.}

\noindent Practically it means that the saddle-point of the AL optimization
can be used in order to obtain the solutions of the considered optimization
problems.

The convergence properties for the analysis and synthesis-based algorithms
are formulated in the following proposition.

\begin{proposition}
\textit{\label{prop:An_Syn_convergence}}

\textit{(a) If there exists a saddle point }$\mathbf{(y}^{\ast }\mathbf{%
,\omega }^{\ast }\mathbf{,\lambda }^{\ast }\mathbf{)}$\textit{\ of }$L_{%
\text{a}}\left( \mathbf{y,\omega ,\lambda }\right) $\textit{\ (\ref{A_AL_1}%
), then }$\mathbf{y}_{t}\mathbf{\rightarrow y}^{\ast }\mathbf{,\omega }_{t}%
\mathbf{\rightarrow \omega }^{\ast }\mathbf{,\lambda }_{t}\mathbf{%
\rightarrow \lambda }^{\ast }$\textit{.}

\textit{(b) If there exists a saddle point }$\mathbf{(y}^{\ast }\mathbf{%
,\omega }^{\ast }\mathbf{,u}^{\ast }\mathbf{,\lambda }^{\ast }\mathbf{)}$%
\textit{\ of }$L_{\text{s}}\left( \mathbf{y,\omega ,u},\mathbf{\lambda }%
\right) $\textit{\ (\ref{S_AL_11}), then }$\mathbf{y}_{t}\mathbf{\rightarrow
y}^{\ast }\mathbf{,\omega }_{t}\mathbf{\rightarrow \omega }^{\ast }\mathbf{,u%
}_{t}\mathbf{\rightarrow u}^{\ast }\mathbf{,\lambda }_{t}\mathbf{\rightarrow
\lambda }^{\ast }$\textit{.}

\textit{On the other hand, if no such saddle point exists, then at least one
of the sequences }$\mathbf{\{y}_{t}\mathbf{\}}\ $\textit{or }$\mathbf{%
\{\lambda }_{t}\mathbf{\}}$\textit{\ must be unbounded.}
\end{proposition}

The proof is given in Appendix \ref{app:convergence}.

\subsection{IDD-BM3D algorithm}

\begin{proposition}
\label{prop:Decouple_convergence}\textit{For any set of parameters }$\sigma
,\tau ,\gamma ,\xi $\textit{\ the sequence }$\left( \mathbf{y}_{t},\mathbf{%
\omega }_{t}\right) $\textit{\ generated by the IDD-BM3D algorithm with
equal group weights }$g_{r}$\textit{, converges to the fixed point }$\left( 
\mathbf{y}^{\ast },\mathbf{\omega }^{\ast }\right) $\textit{\ defined by the
equations (\ref{var_decoupled}), if the fixed point exists.}
\end{proposition}

The proof of the proposition is given in Appendix \ref{app:convergence}. It
is not required that the fixed point is unique. Depending on a starting
point $(\mathbf{y}_{0},\mathbf{\omega }_{0})$ the limit point of the
algorithm can be different but should satisfy the fixed point equations.

\section{Implementation}

\label{sec:Implementations}

\emph{Grouping and frame operators}. To build the groups, we use the
block-matching procedure from \cite{BM3D_TIP} and apply it to the image
reconstructed by the BM3DDEB deblurring algorithm \cite{DEBBM3D_SPIE2008}.
The found locations of the similar blocks constitute the set $J$ that is
necessary to construct the analysis and synthesis frames. Multiplications
against the matrices $\mathbf{\Phi ,\Phi }^{T}\mathbf{,\Psi }$ and $\mathbf{%
\Psi }^{T}$ are calculated efficiently since all of them involve only
groupwise separable 3-D transformations of the data (possibly with some
averaging of the estimates). In our experiments the 3-D transform is
performed by first applying the 2-D discrete sine transform (DST) to each
block in the group followed by the 1-D Haar transform applied along the
third dimension of the group. The image block size is $4\times 4$, and the
number of blocks in the group is $8$.

\emph{Choice of the group weights.} Since image blocks are overlapping, for
each pixel we obtain several estimates. The weighted averaging can be used
to improve the final aggregated estimate. For the one-step (non-iterative)
algorithms the weights can be adaptively selected so to minimize the
variance of the final aggregated estimate, based on the variance of each of
the estimates (e.g. \cite{Guleryuz-OvercompleteDictionaries}, \cite{BM3D_TIP}%
, \cite{DEBBM3D_SPIE2008}{\small )}. In the considered iterative algorithms
the influence of the weights on the final estimate is complex, and deriving
a formula for the optimal weights is rather involved. Instead, following the
idea of the sparse representations, we suggest giving the preference to the
estimates obtained from the sparser groups. In our implementations we use
weights inversely proportional to the number of significant spectrum
coefficients of the groups $g_{r}=1/\left\Vert \mathfrak{Th}_{\epsilon
}\left( \mathbf{\omega }_{r}\right) \right\Vert _{0}$, where significant
coefficients are found by the hard thresholding of the group spectra using a
small threshold $\epsilon $.

The grouping and the adaptive group weights are calculated only once, using
the initial image estimate $\mathbf{y}_{\text{init}}$ and remain unchanged
through the subsequent iterations.

\emph{Choice of the regularization parameters}. The parameters $\tau ,\gamma
,\xi $ are optimized to provide best reconstruction quality. Optimization
has been performed separately for each algorithm and each deblurring
scenario. The parameter $\beta $ is always set to 1.

\emph{Initialization.} We experimentally confirmed the convergence to an
asymptotic solution that is independent of the initialization $\mathbf{y}%
_{0} $ and $\mathbf{\omega }_{0}$. Nevertheless, initialization with a
better estimate, for example with the reconstruction obtained by BM3DDEB
(which we also use to define grouping) results in a much faster convergence.

\emph{Solution of the large-scale linear equations}. All proposed algorithms
contain steps involving solution of large-scale linear equations. For a
circular shift-invariant blur operator, the solution of the equations\textbf{%
\ }(\ref{S_AL_3}\textbf{) }and\textbf{\ (}\ref{eq:Dec_AL_1}\textbf{) }can be
calculated in the Fourier domain using the FFT. The more complex equations (%
\ref{A-AL_3}) and (\ref{S_AL_3_2}) are solved using the conjugate gradient
method. The conjugate gradient method allows avoiding explicit calculations
of the matrices $\Phi ^{T}\Phi $\ and $\Psi ^{T}\Psi $, since it requires
only evaluating products of these matrices against vectors.

\emph{Practical considerations. }The two steps of the IDD-BM3D algorithm can
be merged into a single one%
\begin{equation*}
\mathbf{y}_{t+1}=\mathcal{F}^{-1}\left( \dfrac{\mathcal{F}^{\ast }\left( 
\mathbf{h}\right) \circ \mathcal{F}\left( \mathbf{z}\right) +\frac{\sigma
^{2}}{\gamma }\mathcal{F}\left( \mathbf{\Psi }\mathfrak{Th}_{\tau \xi
}\left( \mathbf{\Phi }\mathbf{y}_{t}\right) \right) }{\left\vert \mathcal{F}%
\left( \mathbf{h}\right) \right\vert ^{2}+\frac{\sigma ^{2}}{\gamma }}%
\right) ,
\end{equation*}%
where the analysis-thresholding-synthesis operation $\Psi \mathfrak{Th}%
_{\tau \xi }\left( \mathbf{\Phi }\mathbf{y}_{t}\right) $\ can be calculated
groupwise without need to obtain the whole spectrum $\omega _{t}$\
explicitly. Here $\mathbf{h}$\ denotes the vectorized blurring kernel
corresponding to the blur operator $\mathbf{A}$, and '$\circ $' stands for
the elementwise product of two vectors. The operator $\mathcal{F}\left(
\cdot \right) $ reshapes the input vector into a 2-D array, performs 2-D FFT
and vectorizes the obtained result. $\mathcal{F}^{-1}\left( \cdot \right) $
works analogously, performing inverse FFT.

\emph{Complexity. }Application of the frame operators is the most
computationally expensive part of the proposed algorithms. However, due to
their specific structure, the complexity of the frame operators $\mathbf{%
\Phi }$ and $\mathbf{\Psi }$ is growing only linearly with respect to the
number of the pixels in the image. To give an estimate of the complexity of
the IDD-BM3D algorithm, we mention that, on a $256\times 256$ image, one
iteration takes about 0.35 seconds, and about 50 iterations are typically
sufficient. This timing has been done on dual core 2.6 GHz processor for an
implementation where the computationally most intensive parts have been
written in C++.

\section{Experiments}

\label{sec:Experiments}

We consider six deblurring scenarios used as the benchmarks in many
publications (e.g., \cite{PORTILLA_2007} and \cite{DEBBM3D_SPIE2008}). The
blur point spread function (PSF) $h\left( x_{1},x_{2}\right) $ and the
variance of the noise $\sigma ^{2}$ for each scenario are summarized in
Table \ref{tbl_scenarios_parameters}. PSFs are normalized so that $\sum h=1$%
. Each of the scenarios was tested with the four standard images:\ \emph{%
Cameraman, Lena, House} and \emph{Barbara}.

%TCIMACRO{\TeXButton{B}{\begin{table}[t] \centering}}%
%BeginExpansion
\begin{table}[t] \centering%
%EndExpansion
\begin{tabular}{|c|c|c|}
\hline
Scenario & PSF & $\sigma ^{2}$ \\ \hline
$1$ & \multicolumn{1}{|l|}{$1/\left( 1+x_{1}^{2}+x_{2}^{2}\right) ,$ $%
x_{1},x_{2}=-7,\dots ,7$} & $2$ \\ \hline
$2$ & \multicolumn{1}{|l|}{$1/\left( 1+x_{1}^{2}+x_{2}^{2}\right) ,$ $%
x_{1},x_{2}=-7,\dots ,7$} & $8$ \\ \hline
$3$ & \multicolumn{1}{|l|}{$9\times 9$ uniform} & $\approx 0.3$ \\ \hline
$4$ & \multicolumn{1}{|l|}{$\left[ 1\text{ }4\text{ }6\text{ }4\text{ }1%
\right] ^{T}\left[ 1\text{ }4\text{ }6\text{ }4\text{ }1\right] /256$} & $49$
\\ \hline
$5$ & \multicolumn{1}{|l|}{Gaussian with $std=1.6$} & $4$ \\ \hline
$6$ & \multicolumn{1}{|l|}{Gaussian with $std=0.4$} & $64$ \\ \hline
\end{tabular}%
\caption{Blur PSF and noise variance used in each scenario.}\label%
{tbl_scenarios_parameters}%
%TCIMACRO{\TeXButton{E}{\end{table}}}%
%BeginExpansion
\end{table}%
%EndExpansion

%TCIMACRO{\TeXButton{B}{\begin{figure}[t!]\centering}}%
%BeginExpansion
\begin{figure}[t!]\centering%
%EndExpansion
\includegraphics[height=2.4951in,width=3.116in]{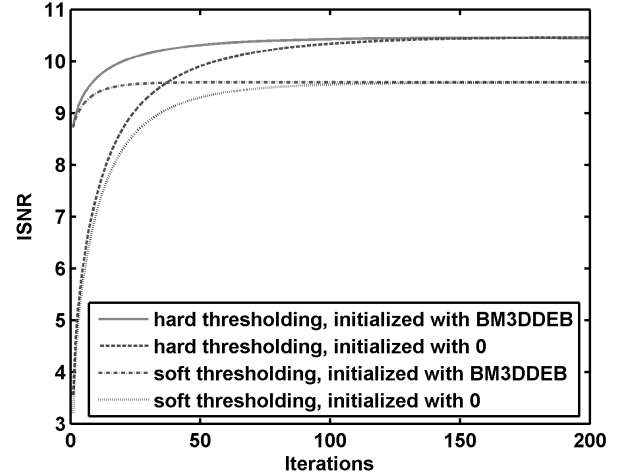}%
\caption{Change of the
ISNR with iterations for the different setups of the IDD-BM3D algorithm. Deblurring
of \textit{Cameraman} image, scenario 3. }\label{fig:Convergence}%
%TCIMACRO{\TeXButton{E}{\end{figure}}}%
%BeginExpansion
\end{figure}%
%EndExpansion

\subsection{Experiment 1 - comparison of the proposed algorithms}

All three proposed algorithms, namely: analysis-based, synthesis-based and
IDD-BM3D are evaluated in the scheme with the soft thresholding and unit
group weights ($g_{r}=1$). Additionally, the IDD-BM3D algorithm is tested
with the adaptive group weights ($g_{r}=1/\left\Vert \mathfrak{Th}_{\epsilon
}\left( \mathbf{\omega }_{r}\right) \right\Vert _{0}$) using the soft and
hard thresholdings.

In Table \ref{tbl_ISNR_exp1} we present improvement of signal-to-noise ratio
(ISNR) values achieved by each algorithm for the \emph{Cameraman }image.
From these values we can conclude that the synthesis-based algorithm
performs essentially worse than the IDD-BM3D algorithm, with the
analysis-based algorithm being in-between. We can also see that the adaptive
weights indeed provide a noticeable restoration improvement. Finally,
comparing the last two rows, we conclude that hard thresholding enables
better results than the soft thresholding, and combined with the adaptive
weights it provides the best results among the considered algorithms.

Convergence properties of the IDD-BM3D algorithm are demonstrated in Figure %
\ref{fig:Convergence}.

The experiments with the IDD-BM3D algorithm can be reproduced using the
Matlab program available as a part of the BM3D package\footnote{%
http://www.cs.tut.fi/\symbol{126}foi/GCF-BM3D}.

\subsection{Experiment 2 - comparison with the state of the art}

Table \ref{tbl_ISNR_exp2} presents a comparison of the IDD-BM3D algorithm
versus a number of algorithms including the current state of the art.\ The
ISNR values for ForWaRD \cite{FORWARD_2004}, SV-GSM \cite{PORTILLA_2007},
SA-DCT \cite{IMDEBLUR_SADCT06} and BM3DDEB \cite{DEBBM3D_SPIE2008} are taken
from our previous paper \cite{DEBBM3D_SPIE2008}, while the results for
L0-AbS \cite{Portilla-ICIP-2009}, TVMM \cite{Dias-2009-TVMM}, CGMK \cite%
{Chantas2010} are obtained by the software available online. We use the
default parameters suggested by the authors of the algorithms. The IDD-BM3D
algorithm in this comparison employs the hard thresholding and the adaptive
weights.

The proposed IDD-BM3D algorithm\ provides the best results with significant
advantage over closest competitors. Particularly interesting is the
comparison against the BM3DDEB algorithm. BM3DDEB is a two-stage
non-iterative algorithm. On the first stage it utilizes the BM3D image
modeling to obtain the initial estimate, which is then used on the second
stage for an empirical Wiener filtering. Better performance of the IDD-BM3D
algorithm demonstrates that considered decoupled formulation (\ref%
{var_decoupled}) enables more effective exploiting of the BM3D-modeling than
the two-stage approach of BM3DDEB.

\begin{figure*}[pt!]
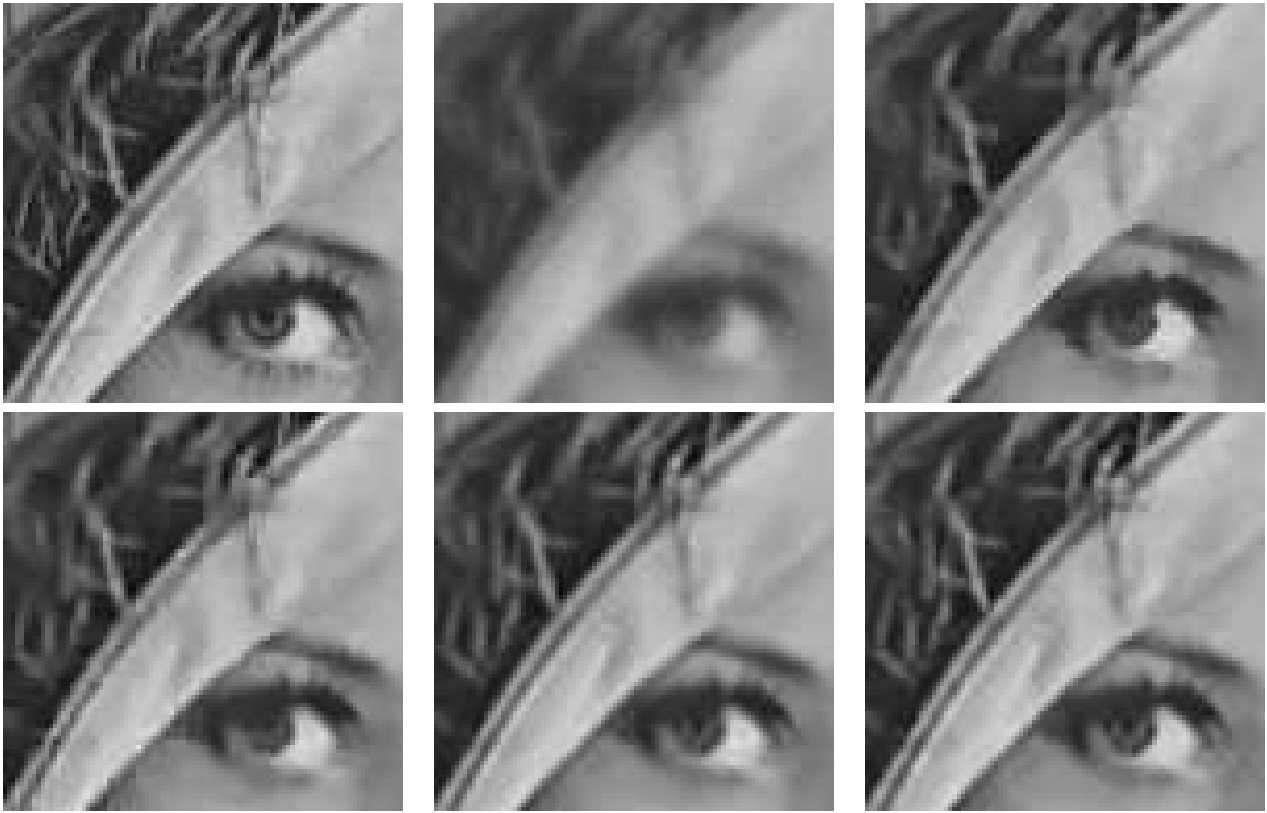
\centering%
\begin{tabular}{lll}
	\includegraphics[height=2.088in,width=2.088in]{true_exp2_lena_crop.eps} & 
	\includegraphics[height=2.088in,width=2.088in]{blurry_exp2_lena_crop.eps} &
	\includegraphics[height=2.088in,width=2.088in]{cgmk_exp2_lena_crop.eps} \\ 
	\includegraphics[height=2.088in,width=2.088in]{l0-abs_exp2_lena_crop.eps} &
	\includegraphics[height=2.088in,width=2.088in]{bm3ddeb_exp2_lena_crop.eps} &
	\includegraphics[height=2.088in,width=2.088in]{iddbm3d_exp2_lena_h_rwi_rwi_fwfg_0__200_crop.eps}
\end{tabular}%
\caption{Deblurring of the \textit{Lena} image, scenario 2. From left to right and from top to bottom are
presented zoomed fragments of the following images: original, blurred noisy,
reconstructed by CGMK \cite{Chantas2010} (ISNR 5.37), L0-AbS \cite{Portilla-ICIP-2009} (ISNR 5.71), DEB-BM3D \cite{DEBBM3D_SPIE2008} (ISNR 6.53) and by proposed IDD-BM3D method (ISNR 6.61).}%
\label{fig_Lena}%
\end{figure*}

The visual quality of some of the restored images can be evaluated from
Figures \ref{fig_Cameraman} and \ref{fig_Lena}, where for a comparison we
show results by the closest competitors \cite{Chantas2010}, \cite%
{Portilla-ICIP-2009} and \cite{DEBBM3D_SPIE2008}. One can see that the
proposed algorithm is able to suppress the ringing artifacts better than
BM3DDEB and provides sharper image edges. This latter effect is achieved in
particular due to the smaller block size used in IDD-BM3D compared to
BM3DDEB.

\section{Discussion}

\label{sec:Discussion}

In the experiments of the previous section we observed a clear advantage of
the IDD-BM3D algorithm over the analysis-based one. This result is rather
surprising, since in the case of the tight frames the IDD-BM3D and the
analysis-based algorithms are almost identical.

Indeed, if we assume that $\left\{ \mathbf{\phi }_{n}\right\} $ is a tight
frame and require that all group weights will be equal, then $\mathbf{\Phi }%
^{T}\mathbf{\Phi }=\alpha \mathbf{I}$ and $\mathbf{\Psi }=\left( \mathbf{%
\Phi }^{T}\mathbf{\Phi }\right) ^{-1}\mathbf{\Phi }^{T}=\alpha ^{-1}\mathbf{%
\Phi }^{T}$. Substituting these expressions into equation (\ref{A-AL_3}) of
the analysis-based algorithm we obtain%
\begin{equation*}
\left( \frac{1}{\sigma ^{2}}\mathbf{A}^{T}\mathbf{A}+\frac{\alpha }{\gamma }%
\mathbf{I}\right) \cdot \mathbf{y=}\frac{1}{\sigma ^{2}}\mathbf{\mathbf{A}%
^{T}\mathbf{z}}+\frac{\alpha }{\gamma }\mathbf{\Psi }\left( \mathbf{\omega }+%
\mathbf{\lambda }\right) .
\end{equation*}%
Comparing it with the equation (\ref{eq:Dec_AL_1}) we see that up to the
presence of the Lagrange multipliers the analysis-based algorithm is
identical to the IDD-BM3D algorithm. This observation rises a question: what
makes the algorithms behave differently when the frame is not tight?

To find an answer, let us look again at the equation (\ref{A-AL_3}). Its
solution requires inversion of the matrix $\frac{1}{\sigma ^{2}}\mathbf{A}%
^{T}\mathbf{A}+\frac{1}{\gamma }\mathbf{\Phi }^{T}\mathbf{\Phi }$, whoes
condition number depends not only on the properties of the blur operator but
also on the properties of the frame. In the case of the non-tight analysis
BM3D-frame, $\mathbf{\Phi }^{T}\mathbf{\Phi }$ is a diagonal matrix, its
entries are defined by the data grouping and count number of times each
pixel appears in different groups. Experiments demonstrate that the
variation of these entries can be very large (up to hundreds times). The
large differences in magnitude of the diagonal elements of $\mathbf{\Phi }%
^{T}\mathbf{\Phi }$ make the matrix $\frac{1}{\sigma ^{2}}\mathbf{A}^{T}%
\mathbf{A}+\frac{1}{\gamma }\mathbf{\Phi }^{T}\mathbf{\Phi }$
ill-conditioned and result in degradation of image reconstruction compared
to IDD-BM3D.

Presence of the matrix $\mathbf{\Phi }^{T}\mathbf{\Phi }$ in the
reconstruction formulas is inevitable as long as one uses criterion
containing norms both for the image and spectrum domain. Formulation based
on the Nash equilibrium allows to overcome this problem and have norms only
from one domain in each criterion.

\section{Conclusions}

\label{sec:Conclusions}

The frame based formulation opens new perspectives for the use of BM3D
modeling within the variational reconstruction techniques. The developed
deblurring algorithm demonstrates state-of-the-art performance, confirming a
valuable potential of BM3D-frames as an advanced image modeling tool. For
non-tight frames, we argue the validity of image reconstruction by
minimizing a single objective function and propose an alternative
formulation, based on Nash equilibrium problem.

%TCIMACRO{\TeXButton{B}{\begin{table*}[tp] \centering}}%
%BeginExpansion
\begin{table*}[tp] \centering%
%EndExpansion
\begin{tabular}{|r|c|c|c|c|c|c|c|c|}
\hline
\multicolumn{3}{|c|}{Method} & \multicolumn{6}{|c|}{Scenario} \\ \hline
& Thresh. & Weights $g_{r}$ & 1 & 2 & 3 & 4 & 5 & 6 \\ \hline\hline
\multicolumn{9}{|c|}{Cameraman (256x256)} \\ \hline
\multicolumn{1}{|l|}{BSNR} &  &  & 31.87 & 25.85 & 40.00 & 18.53 & 29.19 & 
17.76 \\ \hline
\multicolumn{1}{|l|}{Input PSNR} &  &  & 22.23 & 22.16 & 20.76 & 24.62 & 
23.36 & 29.82 \\ \hline\hline
\multicolumn{1}{|l|}{Synthesis} & soft & unit & 6.30 & 4.60 & 7.88 & 2.06 & 
2.98 & 2.84 \\ \hline
\multicolumn{1}{|l|}{Analysis} & soft & unit & 7.88 & 5.75 & 9.22 & 3.00 & 
3.67 & 3.92 \\ \hline
\multicolumn{1}{|l|}{IDD-BM3D} & soft & unit & 8.17 & 6.17 & 9.38 & 3.17 & 
3.83 & 4.12 \\ \hline
\multicolumn{1}{|l|}{IDD-BM3D} & soft & adaptive & 8.41 & 6.41 & 9.59 & 3.38
& 3.98 & 4.14 \\ \hline
\multicolumn{1}{|l|}{IDD-BM3D} & hard & adaptive & \textbf{8.85} & \textbf{%
7.12} & \textbf{10.45} & \textbf{3.98} & \textbf{4.31} & \textbf{4.89} \\ 
\hline
\end{tabular}%
\caption{Comparison of the output ISNR [dB] of the proposed deblurring algorithms. Row
corresponding to ``Input PSNR'' contain PSNR [dB] of the input blurry
images). Blurred signal-to-noise ratio (BSNR) is defined as $10log_{10}\left(var\left(\bold{Ay}\right)/N\sigma^2\right)$, where $var()$ is the variance.}%
\label{tbl_ISNR_exp1}%
%TCIMACRO{\TeXButton{E}{\end{table*}}}%
%BeginExpansion
\end{table*}%
%EndExpansion

%TCIMACRO{\TeXButton{B}{\begin{table*}[tp] \centering}}%
%BeginExpansion
\begin{table*}[tp] \centering%
%EndExpansion
\begin{tabular}{|c||c|c|c|c|c|c||c|c|c|c|c|c|}
\hline
& \multicolumn{6}{||c}{Scenario} & \multicolumn{6}{||c|}{Scenario} \\ \hline
\multicolumn{1}{|r||}{} & 1 & 2 & 3 & 4 & 5 & 6 & 1 & 2 & 3 & 4 & 5 & 6 \\ 
\hline\hline
Method & \multicolumn{6}{||c||}{Cameraman (256x256)} & \multicolumn{6}{||c|}{
House (256x256)} \\ \hline\hline
\multicolumn{1}{|l||}{BSNR$\hspace{-0.3cm}$} & 31.87 & 25.85 & 40.00 & 18.53
& 29.19 & 17.76 & 29.16 & 23.14 & 40.00 & 15.99 & 26.61 & 15.15 \\ \hline
\multicolumn{1}{|l||}{Input PSNR} & 22.23 & 22.16 & 20.76 & 24.62 & 23.36 & 
29.82 & 25.61 & 25.46 & 24.11 & 28.06 & 27.81 & 29.98 \\ \hline\hline
\multicolumn{1}{|l||}{ForWaRD \cite{FORWARD_2004}} & 6.76 & 5.08 & 7.34 & 
2.40 & 3.14 & 3.92 & 7.35 & 6.03 & 9.56 & 3.19 & 3.85 & 5.52 \\ \hline
\multicolumn{1}{|l||}{SV-GSM \cite{PORTILLA_2007}} & 7.45 & 5.55 & 7.33 & 
2.73 & 3.25 & 4.19 & 8.64 & 7.03 & 9.04 & 4.30 & 4.11 & 6.02 \\ \hline
\multicolumn{1}{|l||}{SA-DCT \cite{IMDEBLUR_SADCT06}} & 8.11 & 6.33 & 8.55 & 
\emph{3.37} & 3.72 & \emph{4.71} & 9.02 & 7.74 & 10.50 & 4.99 & 4.65 & 5.96
\\ \hline
\multicolumn{1}{|l||}{BM3DDEB \cite{DEBBM3D_SPIE2008}} & \emph{8.19} & \emph{%
6.40} & 8.34 & 3.34 & \emph{3.73} & 4.70 & \emph{9.32} & \emph{8.14} & \emph{%
10.85} & \emph{5.13} & 4.56 & \textbf{7.21} \\ \hline
\multicolumn{1}{|l||}{L0-AbS \cite{Portilla-ICIP-2009}} & 7.70 & 5.55 & 9.10
& 2.93 & 3.49 & 1.77 & 8.40 & 7.12 & 11.06 & 4.55 & 4.80 & 2.15 \\ \hline
\multicolumn{1}{|l||}{TVMM \cite{Dias-2009-TVMM}} & 7.41 & 5.17 & 8.54 & 2.57
& 3.36 & 1.30 & 7.98 & 6.57 & 10.39 & 4.12 & 4.54 & 2.44 \\ \hline
\multicolumn{1}{|l||}{CGMK \cite{Chantas2010}} & 7.80 & 5.49 & \emph{9.15} & 
2.80 & 3.54 & 3.33 & 8.31 & 6.97 & 10.75 & 4.48 & \emph{4.97} & 4.59 \\ 
\hline
\multicolumn{1}{|l||}{IDD-BM3D} & \textbf{8.85} & \textbf{7.12} & \textbf{%
10.45} & \textbf{3.98} & \textbf{4.31} & \textbf{4.89} & \textbf{9.95} & 
\textbf{8.55} & \textbf{12.89} & \textbf{5.79} & \textbf{5.74} & \emph{7.13}
\\ \hline\hline
& \multicolumn{6}{||c||}{Lena (512x512)} & \multicolumn{6}{||c|}{Barbara
(512x512)} \\ \hline\hline
\multicolumn{1}{|l||}{BSNR} & 29.89 & 23.87 & 40.00 & 16.47 & 27.18 & 15.52
& 30.81 & 24.79 & 40.00 & 17.35 & 28.07 & 16.59 \\ \hline
\multicolumn{1}{|l||}{Input PSNR} & 27.25 & 27.04 & 25.84 & 28.81 & 29.16 & 
30.03 & 23.34 & 23.25 & 22.49 & 24.22 & 23.77 & 29.78 \\ \hline\hline
\multicolumn{1}{|l||}{ForWaRD \cite{FORWARD_2004}} & 6.05 & 4.90 & 6.97 & 
2.93 & 3.50 & 5.42 & 3.69 & 1.87 & 4.02 & 0.94 & 0.98 & 3.15 \\ \hline
\multicolumn{1}{|l||}{SV-GSM \cite{PORTILLA_2007}} & - & - & - & - & - & - & 
6.85 & 3.80 & 5.07 & \textbf{1.94} & \textbf{1.36} & 5.27 \\ \hline
\multicolumn{1}{|l||}{SA-DCT \cite{IMDEBLUR_SADCT06}} & 7.55 & 6.10 & 7.79 & 
4.49 & 4.08 & 5.84 & 5.45 & 2.54 & 4.79 & 1.31 & 1.02 & 3.83 \\ \hline
\multicolumn{1}{|l||}{BM3DDEB \cite{DEBBM3D_SPIE2008}} & \emph{7.95} & \emph{%
6.53} & \emph{7.97} & \emph{4.81} & \emph{4.37} & \textbf{6.40} & \textbf{%
7.80} & \emph{3.94} & \emph{5.86} & \emph{1.90} & \emph{1.28} & \textbf{5.80}
\\ \hline
\multicolumn{1}{|l||}{L0-AbS \cite{Portilla-ICIP-2009}} & 6.66 & 5.71 & 7.79
& 4.09 & 4.22 & 1.93 & 3.51 & 1.53 & 3.98 & 0.73 & 0.81 & 1.17 \\ \hline
\multicolumn{1}{|l||}{TVMM \cite{Dias-2009-TVMM}} & 6.36 & 4.98 & 7.47 & 3.52
& 3.61 & 2.79 & 3.10 & 1.33 & 3.49 & 0.41 & 0.75 & 0.59 \\ \hline
\multicolumn{1}{|l||}{CGMK \cite{Chantas2010}} & 6.76 & 5.37 & 7.86 & 3.49 & 
3.93 & 4.46 & 2.45 & 1.34 & 3.55 & 0.44 & 0.81 & 0.38 \\ \hline
\multicolumn{1}{|l||}{IDD-BM3D} & \textbf{7.97} & \textbf{6.61} & \textbf{%
8.91} & \textbf{4.97} & \textbf{4.85} & \emph{6.34} & \emph{7.64} & \textbf{%
3.96} & \textbf{6.05} & 1.88 & 1.16 & \emph{5.45} \\ \hline
\end{tabular}%
\caption{Comparison of the output ISNR [dB] of deconvolution methods (row
corresponding to ``Input PSNR'' contain PSNR [dB] of the input blurry
images).}\label{tbl_ISNR_exp2}%
%TCIMACRO{\TeXButton{E}{\end{table*}}}%
%BeginExpansion
\end{table*}%
%EndExpansion

\begin{figure*}[pt!]\centering%
\begin{tabular}{lll}

	\includegraphics[height=2.088in,width=2.088in]{true_exp3_cameraman_crop.eps} & 
	\includegraphics[height=2.088in,width=2.088in]{blurry_exp3_cameraman_crop.eps} & 
	\includegraphics[height=2.088in,width=2.088in]{cgmk_exp3_cameraman_crop.eps} \\ 
	\includegraphics[height=2.088in,width=2.088in]{l0-abs_exp3_cameraman_crop.eps} & 
	\includegraphics[height=2.088in,width=2.088in]{bm3ddeb_exp3_cameraman_crop.eps} & 
	\includegraphics[height=2.088in,width=2.088in]{iddbm3d_exp3_cameraman_h_rwi_rwi_fwfg_0__200_crop.eps}
\end{tabular}%
\caption{Deblurring of the\textit{Cameraman} image, scenario 3. From left to right and from top to bottom are
presented zoomed fragments of the following images: original, blurred noisy, reconstructed by CGMK \cite{Chantas2010} (ISNR 9.15), L0-AbS \cite{Portilla-ICIP-2009} (ISNR 9.10), 
DEB-BM3D \cite{DEBBM3D_SPIE2008} (ISNR 8.34) and by proposed IDD-BM3D method (ISNR 10.45).}%
\label{fig_Cameraman}%
\end{figure*}%

\appendices

%TCIMACRO{\TeXButton{Appendix A}{\section{}}}%
%BeginExpansion
\section{}%
%EndExpansion

\label{app:FrameProps}

\subsection{Proof of Proposition \protect\ref{ref_propFrames}}

The proof is based on use of the following Kronecker matrix product formulas.

If $\mathbf{A}$ is an $m\times n$ matrix and $\mathbf{B}$ is a $p\times q$
matrix, then the Kronecker product $\mathbf{A}\otimes \mathbf{B}$ is the $%
mp\times nq$ block matrix and 
\begin{eqnarray*}
(\mathbf{A}\otimes \mathbf{B})(\mathbf{C}\otimes \mathbf{D}) &=&\mathbf{AC}%
\otimes \mathbf{BD}, \\
(\mathbf{A}\otimes \mathbf{B})^{T} &=&\mathbf{A}^{T}\otimes \mathbf{B}^{T},
\\
(\mathbf{A}\otimes \mathbf{B})^{-1} &=&\mathbf{A}^{-1}\otimes \mathbf{B}%
^{-1}.
\end{eqnarray*}%
Also, matrix equation $\mathbf{AXB}=\mathbf{C}$ can be vectorized columnwise
with respect to $\mathbf{X}$ and $\mathbf{C}$ as following 
\begin{equation*}
(\mathbf{B}^{T}\otimes \mathbf{A})vect\left( \mathbf{X}\right) =vect\left( 
\mathbf{C}\right) .
\end{equation*}

To simplify notation we denote $\mathbf{G}=\left( \mathbf{D}_{1}\otimes 
\mathbf{D}_{1}\right) $. Then the formula (\ref{Prop_1_1}) from Proposition %
\ref{ref_propFrames} is proved as following 
\begin{eqnarray*}
&&\mathbf{\Phi }^{T}\mathbf{\Phi }=\sum_{r}\mathbf{\Phi }_{r}^{T}\mathbf{%
\Phi }_{r}= \\
&&\sum_{r}\sum_{j\in J_{r}}\sum_{j^{\prime }\in J_{r}}(\mathbf{d}%
_{j}^{T}\otimes \mathbf{P}_{j}^{T}\mathbf{G}^{T})(\mathbf{d}_{,j^{\prime
}}\otimes \mathbf{GP}_{j^{\prime }})= \\
&&\sum_{r}\sum_{j\in J_{r}}\sum_{j^{\prime }\in J_{r}}(\mathbf{d}_{j}^{T}%
\mathbf{d}_{j^{\prime }})\otimes (\mathbf{P}_{j}^{T}\mathbf{G}^{T}\mathbf{GP}%
_{j})= \\
&&\sum_{r}\sum_{j\in J_{r}}\sum_{j^{\prime }\in J_{r}}\delta _{j,j^{\prime }}%
\mathbf{P}_{j}^{T}\cdot \mathbf{I}\cdot \mathbf{P}_{j^{\prime
}}=\sum_{r}\sum_{j\in J_{r}}\mathbf{P}_{j}^{T}\mathbf{P}_{j}.
\end{eqnarray*}

Proof of the formula (\ref{Prop_1_1_1}):%
\begin{eqnarray*}
&&\mathbf{\Psi \Psi }^{T}=\left( \mathbf{W}^{-1}\cdot \lbrack g_{1}\mathbf{%
\Psi }_{1},\ldots ,g_{R}\mathbf{\Psi }_{R}]\right) \times \\
&&\qquad \qquad \qquad \left( \mathbf{W}^{-1}\cdot \lbrack g_{1}\mathbf{\Psi 
}_{1},\ldots ,g_{R}\mathbf{\Psi }_{R}]\right) ^{T}= \\
&&\mathbf{W}^{-1}[g_{1}\mathbf{\Psi }_{1},\ldots ,g_{R}\mathbf{\Psi }%
_{R}]\cdot \lbrack g_{1}\mathbf{\Psi }_{1},\ldots ,g_{R}\mathbf{\Psi }%
_{R}]^{T}\mathbf{W}^{-1}= \\
&&\mathbf{W}^{-1}\!\!\sum_{r}g_{r}^{2}\!\!\sum_{j\in J_{r}}\sum_{j^{\prime
}\in J_{r}}\!\!\left( \!\mathbf{d}_{j}^{T}\mathbf{\otimes P}_{j}^{T}\mathbf{G%
}^{T}\right) \!\!\left( \mathbf{d}_{j^{\prime }}\mathbf{\otimes GP}%
_{j^{\prime }}\!\right) \!\mathbf{W}^{-1}\!\!= \\
&&\mathbf{W}^{-1}\!\!\sum_{r}g_{r}^{2}\!\!\sum_{j\in J_{r}}\sum_{j^{\prime
}\in J_{r}}\!\!\left( \mathbf{d}_{j}^{T}\mathbf{d}_{j^{\prime }}\right)
\!\otimes \!\left( \mathbf{P}_{j}^{T}\mathbf{G}^{T}\mathbf{GP}_{j^{\prime
}}\right) \!\mathbf{W}^{-1}\!\!= \\
&&\mathbf{W}^{-1}\!\!\sum_{r}g_{r}^{2}\!\!\sum_{j\in J_{r}}\sum_{j^{\prime
}\in J_{r}}\!\!\delta _{j,j^{\prime }}\otimes \left( \mathbf{P}_{j}^{T}%
\mathbf{P}_{j}\right) \!\mathbf{W}^{-1}= \\
&&\mathbf{W}^{-1}\!\!\sum_{r}g_{r}^{2}\!\!\sum_{j\in J_{r}}\mathbf{P}_{j}^{T}%
\mathbf{P}_{j}\mathbf{W}^{-1}= \\
&&\mathbf{W}^{-2}\!\!\sum_{r}g_{r}^{2}\!\!\sum_{j\in J_{r}}\mathbf{P}_{j}^{T}%
\mathbf{P}_{j}\mathbf{.}
\end{eqnarray*}%
The last identity holds since $\sum_{r}g_{r}^{2}\sum_{j\in J_{r}}\mathbf{P}%
_{j}^{T}\mathbf{P}_{j}$ and $\mathbf{W}^{-1}$ are diagonal matrices.

The formula (\ref{Prop_1_2}) in Proposition \ref{ref_propFrames} is valid
since%
\begin{eqnarray*}
&&\mathbf{\Psi \Phi }=\left( \mathbf{W}^{-1}\cdot \lbrack g_{1}\mathbf{\Psi }%
_{1},\ldots ,g_{R}\mathbf{\Psi }_{R}]\right) \times \left[ 
\begin{array}{c}
\mathbf{\Phi }_{1} \\ 
\vdots \\ 
\mathbf{\Phi }_{R}%
\end{array}%
\right] = \\
&&\mathbf{W}^{-1}\sum_{r}\left( g_{r}\mathbf{\Psi }_{r}\right) \mathbf{\Phi }%
_{r}= \\
&&\mathbf{W}^{-1}\sum_{r}\!\left( \!g_{r}\sum_{j\in J_{r}}\mathbf{d}%
_{j}^{T}\otimes \mathbf{P}_{j}^{T}\mathbf{G}^{T}\right) \!\!\left(
\sum_{j^{\prime }\in J_{r}}\mathbf{d}_{j}\otimes \mathbf{GP}_{j^{\prime
}}\!\right) \!\!= \\
&&\mathbf{W}^{-1}\sum_{r}g_{r}\sum_{j\in J_{r}}\sum_{j^{\prime }\in
J_{r}}\left( \mathbf{d}_{j}^{T}\otimes \mathbf{P}_{j}^{T}\mathbf{G}%
^{T}\right) \left( \mathbf{d}_{j}\otimes \mathbf{GP}_{j^{\prime }}\right) =
\\
&&\mathbf{W}^{-1}\sum_{r}g_{r}\sum_{j\in J_{r}}\sum_{j^{\prime }\in
J_{r}}\left( \mathbf{d}_{j}^{T}\mathbf{d}_{j}\right) \otimes \left( \mathbf{P%
}_{j}^{T}\mathbf{G}^{T}\mathbf{GP}_{j^{\prime }}\right) = \\
&&\mathbf{W}^{-1}\sum_{r}g_{r}\sum_{j\in J_{r}}\sum_{j^{\prime }\in
J_{r}}\delta _{ij}\left( \mathbf{P}_{j}^{T}\mathbf{P}_{j^{\prime }}\right) =
\\
&&\mathbf{W}^{-1}\sum_{r}g_{r}\sum_{j\in J_{r}}\mathbf{P}_{j}^{T}\mathbf{P}%
_{j}=I_{N\times N}.
\end{eqnarray*}

%TCIMACRO{\TeXButton{Appendix B}{\section{}}}%
%BeginExpansion
\section{}%
%EndExpansion

\label{app:convergence}

\subsection{Proof of Proposition \protect\ref{prop:An_Syn_convergence}}

Let us consider constrained optimization problem given in the following
general form%
\begin{equation}
\min_{\mathbf{u,v}}\{f\left( \mathbf{u}\right)
+\dsum\limits_{j=1}^{q}g_{j}\left( \mathbf{v}_{j}\right) |\mathbf{Cv}+%
\mathbf{Du=b\}},  \label{CB_4}
\end{equation}%
where $\mathbf{u\in }$\textbf{$%
%TCIMACRO{\U{211d} }%
%BeginExpansion
\mathbb{R}
%EndExpansion
$}$^{m},\mathbf{v=}\left[ \mathbf{v}_{1}^{T},...,\mathbf{v}_{q}^{T}\right]
^{T},\mathbf{v}_{j}\in 
%TCIMACRO{\U{211d} }%
%BeginExpansion
\mathbb{R}
%EndExpansion
^{m_{j}},\mathbf{v}\in 
%TCIMACRO{\U{211d} }%
%BeginExpansion
\mathbb{R}
%EndExpansion
^{\bar{m}},\bar{m}=\sum m_{j},\mathbf{b}\in 
%TCIMACRO{\U{211d} }%
%BeginExpansion
\mathbb{R}
%EndExpansion
^{s},\mathbf{C}$ is of the size $(s\times \bar{m}),\mathbf{D}$ is of the
size $(s\times m)$ and $f(\mathbf{u)}$ is convex. The AL corresponding to
this problem is%
\begin{eqnarray}
L\left( \mathbf{u},\mathbf{v},\mathbf{\lambda }\right) &=&\hspace{-0.2cm}%
f\left( \mathbf{u}\right) +\dsum\limits_{j=1}^{q}g_{j}\left( \mathbf{v}%
_{j}\right) +  \notag \\
&&\hspace{-1cm}\alpha \left\Vert \mathbf{Cv}+\mathbf{Du-b}\right\Vert
_{2}^{2}+\left\langle \mathbf{Cv}+\mathbf{Du-b,\lambda }\right\rangle .
\label{CB_5}
\end{eqnarray}%
The link between the main variable $\mathbf{u}$ and the auxiliary splitting
variable $\mathbf{v}$ is given by the linear equation $\mathbf{Cv}+\mathbf{%
Du=b}$. If $\mathbf{C}$ is the identity matrix, then $\mathbf{v=b-Du}$ and
the convergence of the corresponding iterative algorithm can be obtained
from the Eckstein-Bertsekas's theorem (\cite{Bertsekas-BOOK}, Theorem 8).
However, if $\mathbf{Cv}+\mathbf{Du=b}$ is not resolved with respect to $%
\mathbf{v}$ then the theorem is not applicable in its original form. The
techniques exploited in our paper leads to the relations between the
variables which cannot be resolved with respect to $\mathbf{v}$. In order to
analyze the convergence of the proposed algorithm we use a novel formulation
of the Eckstein-Bertsekas's theorem \cite{Esser} adapted to the general
linear link between the variables $\mathbf{v}$ and $\mathbf{u}$. This new
Eckstein-Bertsekas's theorem is given in the following form \cite{Esser}.

\begin{theorem}
Consider the problem (\ref{CB_4}) where $f$\ and $g_{j}$\ are closed proper
convex functions, $C$\ has full column rank and $f\left( \mathbf{u}\right)
+\left\Vert \mathbf{Du}\right\Vert _{2}^{2}$\ is strictly convex. Let $%
u_{0}\in 
%TCIMACRO{\U{211d} }%
%BeginExpansion
\mathbb{R}
%EndExpansion
^{m},\lambda _{0}\in 
%TCIMACRO{\U{211d} }%
%BeginExpansion
\mathbb{R}
%EndExpansion
^{s}$\ be arbitrary and $\beta >0$. Suppose that there are sequences $%
\left\{ \sigma _{t}^{2}\right\} $\ and $\left\{ \nu _{t}\right\} $\ such
that $\sigma _{t}^{2}\geq 0,\nu _{t}\geq 0$\ and $\sum_{t}\sigma
_{t}^{2}<\infty ,\sum_{t}\nu _{t}<\infty $. Assume that%
\begin{align*}
& \left\Vert \mathbf{v}_{t+1}-\arg \min_{\mathbf{v}}\left\{
\tsum\nolimits_{j=1}^{q}g_{j}\left( \mathbf{v}_{j}\right) +\right. \right. \\
& \qquad \left. \left. +\alpha \left\Vert \mathbf{Cv}+\mathbf{Du}_{t}\mathbf{%
-b}\right\Vert _{2}^{2}+\left\langle \mathbf{Cv,\lambda }_{t}\right\rangle
\right\} \right\Vert _{2}^{2}\leq \sigma _{t}^{2} \\
& \left\Vert \mathbf{u}_{t+1}-\arg \min_{\mathbf{u}}\left\{ %
\vphantom{\left\Vert 1 \right\Vert _{2}^{2}}f(\mathbf{u)}+\right. \right. \\
& \qquad \left. \left. +\alpha \left\Vert \mathbf{Cv}_{t+1}+\mathbf{Du-b}%
\right\Vert _{2}^{2}+\left\langle \mathbf{Du,\lambda }_{t}\right\rangle
\right\} \right\Vert _{2}^{2}\leq \nu _{t}\text{,} \\
& \mathbf{\lambda }_{t+1}=\mathbf{\lambda }_{t}+\beta \left( \mathbf{Cv}%
_{t+1}+\mathbf{Du}_{t+1}-\mathbf{b}\right) .
\end{align*}%
If there exists a saddle point $\left( \mathbf{v}^{\ast }\mathbf{,u}^{\ast }%
\mathbf{,\lambda }^{\ast }\right) $\ for $L\left( \mathbf{u,v,\lambda }%
\right) $\ (\ref{CB_5}), then $v_{t}\rightarrow v^{\ast },u_{t}\rightarrow
u^{\ast },\lambda _{t}\rightarrow \lambda ^{\ast }$. On the other hand, if
no such a saddle point exists, then at least one of the sequences $\left\{ 
\mathbf{u}_{t}\right\} $\ or $\left\{ \mathbf{\lambda }_{t}\right\} $\ must
be unbounded.
\end{theorem}

This formulation of the convergence concerns approximate solutions on each
optimization step, where the parameters $\sigma _{t}^{2}$ and $\nu _{t}$
controls the accuracy at each step. The finite sums $\sum_{t}\sigma
_{t}^{2}<\infty ,\sum_{t}\nu _{t}<\infty $ mean that $\sigma _{t}^{2}\mathbf{%
,}\nu _{t}\rightarrow 0$, i.e. the accuracy should asymptotically improve.

Armed with this theorem we can proceed to the proof of Proposition \ref%
{prop:An_Syn_convergence}.

(a) Comparing the AL (\ref{A_AL_1}) with (\ref{CB_4}) we note that $f\left( 
\mathbf{u}\right) =\dfrac{1}{2\sigma ^{2}}\left\Vert \mathbf{z}-\mathbf{Ay}%
\right\Vert _{2}^{2}$ and the equality $\mathbf{Cv}+\mathbf{Du=b}$ takes the
form $\mathbf{\omega }-\mathbf{\Phi y}=0$, where $\mathbf{\omega }$
corresponds to $\mathbf{v}$ and $\mathbf{u}$ corresponds to $\mathbf{y}$.
Thus, $\mathbf{C=I}_{M\times M}$ and $\mathbf{D=-\Phi }$.

We have two conditions of the theorem to be tested: $\mathbf{C}$ has full
column rank and $f\left( \mathbf{u}\right) \mathbf{+}\left\Vert \mathbf{Du}%
\right\Vert _{2}^{2}$ is strictly convex. In our case, $\mathbf{C=I}%
_{M\times M}$ has full column rank, $\left\Vert \mathbf{Du}\right\Vert
_{2}^{2}=\left\langle \mathbf{\Phi }^{T}\mathbf{\Phi u},\mathbf{u}%
\right\rangle $. Due to (\ref{Prop_1_1}) $\mathbf{\Phi }^{T}\mathbf{\Phi =W>0%
}$, thus $\left\Vert \mathbf{Du}\right\Vert _{2}^{2}$ is strongly convex and
the same holds for $\dfrac{1}{2\sigma ^{2}}\left\Vert \mathbf{z}-\mathbf{Ay}%
\right\Vert _{2}^{2}+\left\Vert \mathbf{Du}\right\Vert _{2}^{2}$. Thus, all
conditions of the theorem are satisfied and the analysis-based algorithm
converges to the saddle-point of the AL (\ref{A_AL_1}), if it exists. It
proves the first part of the proposition.

(b) Comparing the formulation (\ref{S_AL_1}) with (\ref{CB_4}) we note that $%
f\left( \mathbf{u}\right) =\dfrac{1}{2\sigma ^{2}}\left\Vert \mathbf{z}-%
\mathbf{Ay}\right\Vert _{2}^{2}$ and the equality $\mathbf{Cv}+\mathbf{Du=b}$
takes the form $\mathbf{y-\Psi u=0}$ and $\mathbf{\omega -u=0}$. Assuming $%
\mathbf{v}\rightarrow \dbinom{\mathbf{y}}{\mathbf{u}},\mathbf{u}\rightarrow 
\mathbf{\omega }$ these equations give%
\begin{equation*}
\mathbf{C=}\left( 
\begin{array}{cc}
\mathbf{I}_{N\times N} & -\mathbf{\Psi } \\ 
\mathbf{0} & \mathbf{I}_{M\times M}%
\end{array}%
\right) ,\mathbf{D=}\left( 
\begin{array}{c}
0_{N\times M} \\ 
-\mathbf{I}_{M\times M}%
\end{array}%
\right) ,\mathbf{b=0}.
\end{equation*}

The matrix $\mathbf{C}$ is square triangular with elements of the main
diagonal equal to $1$. It has full column rank. For $\left\Vert \mathbf{%
\mathbf{Du}}\right\Vert _{2}^{2}$ we have $\left\Vert \mathbf{Du}\right\Vert
_{2}^{2}\rightarrow \left\Vert \mathbf{\omega }\right\Vert _{2}^{2}$. Thus $%
\left\Vert \mathbf{Du}\right\Vert _{2}^{2}$ is strongly convex and the both
conditions of the theorem are fulfilled. It proves the second part of the
proposition.

\subsection{Proof of Proposition \protect\ref{prop:Decouple_convergence}}

We consider the IDD-BM3D algorithm with soft thresholding and equal group
weights $g_{r}=c,$ $c\in \mathbb{R}^{+},r=1,...,R$. From (\ref{eq:def_PHI}),
(\ref{eq:def_PSI}), (\ref{eq_W}) and (\ref{Prop_1_1}) follows that $\mathbf{%
\Phi }^{T}\mathbf{\Phi =W}$ and $\mathbf{\Psi =W}^{-1}\mathbf{\Phi }^{T}$.

Each iteration of the IDD-BM3D algorithm consists of two steps%
\begin{equation}
\left\{ 
\begin{array}{l}
\mathbf{y}_{t+1}\mathbf{=M}^{-1}\left[ \frac{\gamma }{\sigma ^{2}}\mathbf{%
\mathbf{A}}^{T}\mathbf{z}+\mathbf{\Psi \omega }_{t}\right] , \\ 
\mathbf{\omega }_{t+1}=\mathfrak{Th}_{\tau \xi }\left( \mathbf{\Phi y}%
_{t+1}\right) ,%
\end{array}%
\right.  \label{Decoupl_alg_1}
\end{equation}%
where $\mathbf{M}=\frac{\gamma }{\sigma ^{2}}\mathbf{A}^{T}\mathbf{A}+%
\mathbf{I}>0$.

Introducing the operator $O_{\text{d}}\left( \mathbf{\omega }\right) =%
\mathbf{\Phi M}^{-1}[\frac{\gamma }{\sigma ^{2}}\mathbf{\mathbf{A}}^{T}%
\mathbf{z}+\mathbf{\Psi \omega }]$ and denoting $\mathbf{q}_{t}=\mathbf{\Phi
y}_{t}$ we rewrite (\ref{Decoupl_alg_1}) in a compact form 
\begin{equation}
\left\{ 
\begin{array}{l}
\mathbf{q}_{t+1}=O_{\text{d}}\left( \mathbf{\omega }_{t}\right) , \\ 
\mathbf{\omega }_{t+1}=\mathfrak{Th}_{\tau \xi }\left( \mathbf{q}%
_{t+1}\right) .%
\end{array}%
\right.  \label{Decoupl_alg_2}
\end{equation}

The convergence analysis is based on the technique of nonexpansive
operators. An operator $\mathbf{P}:%
%TCIMACRO{\U{211d} }%
%BeginExpansion
\mathbb{R}
%EndExpansion
^{m}\rightarrow 
%TCIMACRO{\U{211d} }%
%BeginExpansion
\mathbb{R}
%EndExpansion
^{m}$ is called nonexpansive if for any $\mathbf{x},\mathbf{x}^{\prime }\in 
%TCIMACRO{\U{211d} }%
%BeginExpansion
\mathbb{R}
%EndExpansion
^{m}$%
\begin{equation*}
\left\Vert \mathbf{P(\mathbf{x})-P(\mathbf{x}}^{\prime }\mathbf{)}%
\right\Vert _{2}^{2}\leq \left\Vert \mathbf{\mathbf{x}-\mathbf{x}}^{\prime
}\right\Vert _{2}^{2}.
\end{equation*}

It is shown in \cite{Wang-Zhang-2008} (Proposition 3.1) that the soft
thresholding is a nonexpansive operator%
\begin{equation*}
\left\Vert \mathfrak{Th}_{\tau }^{soft}(\mathbf{x})-\mathfrak{Th}_{\tau
}^{soft}\left( \mathbf{x}^{\prime }\right) \right\Vert _{2}^{2}\leq
\left\Vert \mathbf{\mathbf{x}-\mathbf{x}}^{\prime }\right\Vert _{2}^{2},
\end{equation*}%
with equality holding only when 
\begin{equation}
\mathfrak{Th}_{\tau }^{soft}\left( \mathbf{x}\right) -\mathfrak{Th}_{\tau
}^{soft}\left( \mathbf{x}^{\prime }\right) =\mathbf{\mathbf{x}-\mathbf{x}}%
^{\prime }\mathbf{.}  \label{qq}
\end{equation}%
Hence the operator $\mathfrak{Th}_{\tau \xi }\left( \mathbf{\cdot }\right) $
in (\ref{Decoupl_alg_2}) is nonexpansive.

To prove that the operator $O_{\text{d}}$ in (\ref{Decoupl_alg_2}) is also
nonexpansive, we first notice that 
\begin{equation*}
O_{\text{d}}\left( \mathbf{\mathbf{\omega }}\right) -O_{\text{d}}\left( 
\mathbf{\mathbf{\omega }}^{\prime }\right) =\mathbf{\Phi M}^{-1}\mathbf{\Psi 
}\left( \mathbf{\mathbf{\omega }-\mathbf{\omega }}^{\prime }\right) .
\end{equation*}

To find the norm of the matrix $\mathbf{\Phi M}^{-1}\mathbf{\Psi }$ we
evaluate its eigenvalues. For the matrix $\mathbf{\Phi M}^{-1}\mathbf{\Psi }$%
, the corresponding characteristic equation is defined as a determinant of
the equation%
\begin{equation}
(\mathbf{\Phi \mathbf{M}^{-1}W}^{-1}\mathbf{\Phi }^{T}-\lambda \mathbf{I})%
\mathbf{\tilde{v}}=0,  \label{eq2}
\end{equation}%
where $\mathbf{\tilde{v}}\in \mathbb{R}^{M}$ is an eigenvector and $\lambda $
is an eigenvalue. The matrix $\mathbf{\Phi M}^{-1}\mathbf{\Psi }$ has the
size $M\times M$ while its rank is equal to $N$. Thus, $M-N$ eigenvalues of
this matrix are equal to zero. We wish to show that nonzero eigenvalues of $%
\mathbf{\Phi M}^{-1}\mathbf{\Psi }$ coincide with the eigenvalues of the
matrix $\mathbf{M}^{-1}$.

Let us replace in (\ref{eq2}) $\mathbf{\tilde{v}}$ by $\mathbf{\Phi v,v}\in 
\mathbb{R}^{N},$ and multiply the equation (\ref{eq2}) by $\mathbf{W}^{-1}%
\mathbf{\Phi }^{T}$. Then, this equation takes the form%
\begin{equation}
\mathbf{W}^{-1}\mathbf{\Phi }^{T}(\mathbf{\Phi \mathbf{M}^{-1}W}^{-1}\mathbf{%
\Phi }^{T}-\lambda \mathbf{I})\mathbf{\Phi v}=0.  \label{eq3}
\end{equation}%
Multiplication by $\mathbf{W}^{-1}\mathbf{\Phi }^{T}$ in (\ref{eq3}) is
legitimate because it preserves the rank of this system of the linear
equations. Since $\mathbf{W}^{-1}\mathbf{\Phi }^{T}\mathbf{\Phi =I}$, (\ref%
{eq3}) takes the form%
\begin{equation}
(\mathbf{M}^{-1}-\lambda \mathbf{I})\mathbf{v}=0.  \label{eq1}
\end{equation}%
Here $\lambda $ and $\mathbf{v}$ become the eigenvector and eigenvalue for
the matrix $\mathbf{M}^{-1}$. The eigenvalues of the matrix $\mathbf{M}^{-1}=%
\left[ \frac{\gamma }{\sigma ^{2}}\mathbf{A}^{T}\mathbf{A}+\mathbf{I}\right]
^{-1}$ are positive and take values less than or equal to $1$.

The passage from (\ref{eq2}) to (\ref{eq1}) proves that nonzero eigenvalues
of the matrix $\mathbf{\Phi M}^{-1}\mathbf{\Psi }$ are equal to the
eigenvalues of the matrix $\mathbf{M}^{-1}.$Thus all eigenvalues of the
matrix $\mathbf{\Phi M}^{-1}\mathbf{\Psi }$ are nonnegative and take values
less than or equal to 1. Hence, the matrix norm\ $\rho \left( \mathbf{\Phi M}%
^{-1}\mathbf{\Psi }\right) $ is less than or equal to one, and the operator\ 
$O_{\text{d}}$ is nonexpansive due to the inequality%
\begin{multline*}
\left\Vert O_{\text{d}}\left( \mathbf{\omega }\right) -O_{\text{d}}\left( 
\mathbf{\omega }^{\prime }\right) \right\Vert _{2}=\left\Vert \mathbf{\Phi M}%
^{-1}\mathbf{\Psi }\left( \mathbf{\mathbf{\omega }-\mathbf{\omega }}^{\prime
}\right) \right\Vert _{2} \\
\leq \rho \left( \mathbf{\Phi M}^{-1}\mathbf{\Psi }\right) \left\Vert 
\mathbf{\omega }-\mathbf{\omega }^{\prime }\right\Vert _{2}\leq \left\Vert 
\mathbf{\omega }-\mathbf{\omega }^{\prime }\right\Vert _{2}\text{.}
\end{multline*}

Let $\left( \mathbf{y}^{\ast },\mathbf{\omega }^{\ast }\right) $ be a fixed
point of the equations (\ref{Decoupl_alg_1}) and $\Delta \mathbf{y}_{t}%
\mathbf{=y}_{t}-\mathbf{y}^{\ast },$ $\Delta \mathbf{\omega }_{t}=\mathbf{%
\omega }_{t}-\mathbf{\omega }^{\ast },$ $\Delta \mathbf{q}_{t}=\mathbf{\Phi }%
\Delta \mathbf{y}$. Since $\mathfrak{Th}_{\tau \xi }$ and $O_{\text{d}}$ are
nonexpansive operators we have from (\ref{Decoupl_alg_2}) that $||\Delta 
\mathbf{q}_{t+1}||\leq \left\Vert \Delta \mathbf{\omega }_{t}\right\Vert $
and $\left\Vert \Delta \mathbf{\omega }_{t+1}\right\Vert \leq ||\Delta 
\mathbf{q}_{t+1}||$. It follows that $\left\Vert \Delta \mathbf{\omega }%
_{t+1}\right\Vert \leq \left\Vert \Delta \mathbf{\omega }_{t}\right\Vert $
for $\forall t$. Then, the sequence $\mathbf{\omega }_{t+1}$ lies in a
compact region and converging to a limit point, say $\mathbf{\tilde{\omega}}%
, $ $lim_{k\rightarrow \infty }\left\Vert \mathbf{\omega }_{t_{k}}-\mathbf{%
\omega }^{\ast }\right\Vert =||\mathbf{\tilde{\omega}-\mathbf{\omega }^{\ast
}||}$, i.e. a distance from this limit point to a fixed point is bounded. By
the continuity of the operators in (\ref{Decoupl_alg_1}) the same statement
holds for the sequence $\mathbf{y}_{t}$: at least one limit point exists,
denoted as $\mathbf{\tilde{y}}$, and a distance between this limit point and
a fixed point is bounded, $lim_{k\rightarrow \infty }\left\Vert \mathbf{y}%
_{t_{k}}-\mathbf{y}^{\ast }\right\Vert =||\mathbf{\tilde{y}-\mathbf{y}^{\ast
}||}$.

Again due to the continuity of the operators in (\ref{Decoupl_alg_1}) the
limit point is a fixed point. Replacing $\left( \mathbf{y}^{\ast },\mathbf{%
\omega }^{\ast }\right) $ by $\left( \mathbf{\tilde{y}},\mathbf{\tilde{\omega%
}}\right) $ we obtain the convergence of the decoupling algorithm, $%
lim_{k\rightarrow \infty }\left\Vert \mathbf{\omega }_{t_{k}}-\mathbf{\tilde{%
\omega}}\right\Vert =0$ and $lim_{k\rightarrow \infty }\left\Vert \mathbf{y}%
_{t_{k}}-\mathbf{\tilde{y}}\right\Vert =0$. It proves Proposition \ref%
{prop:Decouple_convergence}.

\bibliographystyle{IEEEtran}
\bibliography{deblurring,IEEEabrv}

\end{document}